\documentclass[11pt]{article}
\usepackage{amsfonts, eucal, hyperref}

\newtheorem{Proposition}{Proposition}
\newtheorem{Theorem}{Theorem}
\newtheorem{Lemma}{Lemma}
\newtheorem{Cor}{Corollary}
\newtheorem{Definition}{Definition}
\newtheorem{Remark}{Remark}

\newcommand{\restate}[2]{\medskip \noindent {\bf #1 } {\it #2}}
  
\newcommand{\rmk}{\bigskip \par \noindent {\bf Remark: }}

\newcommand{\exam}[1]{\subsubsection*{Example: #1}}

\newcommand{\e}{\epsilon}

\newcommand{\fg}{\mathfrak{g}}                 
\newcommand{\tfg}{\widetilde {\fg}}                 
\newcommand{\fh}{\mathfrak{h}}                 
                 
\newcommand{\rmd}{\mbox{\rm d}}

\newcommand{\cG}{\mathcal{G}}

\newcommand{\cL}{\mathcal{L}}                 
\newcommand{\cM}{M}

\newcommand{\cS}{\mathcal{S}}                 
                 
\newcommand{\cU}{\mathcal{U}}                 
\newcommand{\cV}{\mathcal{V}}                 
\newcommand{\cW}{\mathcal{W}}

\newcommand{\proj}[1]{\bbP_{\!#1}}
\newcommand{\projh}[1]{\proj{{\rm hor}}}
\newcommand{\prpc}{\proj \parconn}
\newcommand{\prpcf}{\proj \alpha}

\newcommand{\projeh}{\proj{{\fh}}}
\newcommand{\projegh}{\proj{{\Ad g \, \fh}}}

\newcommand{\Ads}[1]{\mathrm{Ad}^*_{#1}} 
\newcommand{\actM}[1]{\widehat \Phi_{#1}} 
\newcommand{\cay}[1]{\mathrm{cay}({#1})} 
\newcommand{\Ann}{\mbox{\rm Ann}\ }

\newcommand{\R}{{\mathbb R}}

\newcommand{\parconn}{\Gamma}
\newcommand{\cod}{\nabla}
\newcommand{\codm}{\nabla_{\!m}}

\newcommand{\wtm}{\widetilde{\mu}}

\newcommand{\wtpp}{\widetilde{\pi}_\phi}                 
\newcommand{\mod}{\qquad \mbox{mod} \quad}
\newcommand{\pff}{partial moving frame}
\newcommand{\pcf}{partial connection form }
\newcommand{\ig}[2]{#1_{\cM}(#2)}
\newcommand{\apply}[1]{\widehat \Phi_{#1}}
\newcommand{\triv}{d^\natural}
\newcommand{\trivl}{d^{\natural_{{\!}_L}}}

\newcommand{\kernel}[1]{\mathrm{ker} \, #1}
\newcommand{\range}[1]{\mathrm{range} \, #1}

\newcommand{\spn}[1]{\mathrm{span} \lcb #1 \rcb}
\newcommand{\bleh}[1]{\rho_\phi(#1)}
\newcommand{\Ad}[1]{\mathrm{Ad}_{#1}}
\newcommand{\ad}[1]{\mathrm{ad}_{#1}}
\newcommand{\ads}[1]{\mathrm{ad}^*_{#1}}
\newcommand{\lp}{\left (}
\newcommand{\rp}{\right )}
\newcommand{\la}{\left \langle}
\newcommand{\ra}{\right \rangle}
\newcommand{\lsb}{\left [}
\newcommand{\rsb}{\right ]}
\newcommand{\lcb}{\left \{}
\newcommand{\rcb}{\right \}}
\newcommand{\sands}{\qquad \mbox{and} \qquad}
\newcommand{\half}{{\textstyle {1 \over 2}}}

\newcommand{\smallfrac}[2]{{\textstyle {#1 \over #2}}}
\newcommand{\dm}{\delta m}
\newcommand{\du}{\delta u}

\newcommand{\for}[1]{(\ref{#1})}

\newcommand{\beq}[1]{\begin{equation}\label{#1}}
\newcommand{\eeq}{\end{equation}}
\newcommand{\beqa}{\begin{eqnarray*}}
\newcommand{\eeqa}{\end{eqnarray*}}
\newcommand{\beqan}{\begin{eqnarray}}
\newcommand{\eeqan}{\end{eqnarray}}
\newcommand{\Prop}[1]{\acapon\begin{Proposition}\label{#1}}
\newcommand{\eProp}{\end{Proposition}}
\newcommand{\Lem}[1]{\acapon\begin{Lemma}\label{#1}}
\newcommand{\eLem}{\end{Lemma}}

\newcommand{\prf}{\noindent {\bf Proof: }}
\newcommand{\prfend}{\ \vrule height6pt width6pt depth0pt \medskip}
\newcommand{\prfof}[1]{\smallskip \noindent {\bf Proof of #1: }}
\newcommand{\rmkend}{ \medskip}

\newcommand{\exampleend}{\ \vrule height6pt width6pt depth0pt \medskip}

\newcommand\rank[1]{\mathrm{rank} \lsb #1 \rsb}

\newcommand{\norm}[1]{\left \Vert #1 \right \Vert}

\newcommand{\bbP}{{\mathbb P}}

\newcommand{\ddeps}{\smallfrac{d\ }{d\e}}

\newcommand{\gmo}{{G_{m_0}}}
\newcommand{\fgmo}{{\fg_{m_0}}}

\newcommand{\interior}[1]{#1 \mathbin{\hbox{\hbox{{\vbox
    {\hrule height.4pt width5pt depth0pt}}}\vrule height5.5pt width.4pt depth0pt}\,}} 
\font\elevenrm=cmr10 at 11pt
\newcommand{\Id}{\hbox{{\rm 1}\kern-3.8pt \elevenrm1}}

\begin{document}

\title{Connections for general group actions}
\author{Debra Lewis\thanks{Mathematics Department, University of California,
Santa Cruz, Santa Cruz, CA 95064, USA. email: lewis@math.ucsc.edu.
Supported in part by NSF Grant DMS 01--04292 and by the UCSC Academic Senate 
Committee on Research},
Nilima Nigam\thanks{Department of Mathematics and Statistics, McGill
University, Montreal, H3A 2K6, Canada. email: nigam@math.mcgill.edu.
Supported in part by NSERC.},
and Peter J.~Olver\thanks{Department of Mathematics, University of Minnesota, 
MN 55455, USA.  
email: olver@math.umn.edu. Supported in part by NSF Grant DMS 01--03944}}

\maketitle

\section{Introduction}

Principal bundles serve as a powerful and elegant geometric framework for
analyzing group actions and symmetry. Beyond their geometric origins, principal
bundles play significant roles in the analysis of mechanical systems with
symmetry, as well as the design of appropriate computational algorithms. A
connection on a principal bundle is defined as an equivariant decomposition of
vectors into infinitesimal group motions and complementary infinitesimal `rigid'
motions, which often facilitates the analysis of a system (see, e.g., 
\cite{KN, steenrod}). In some settings,
the group motions are the crucial information, and often the only information
that is sufficiently simple to lend itself to a rigorous analysis. For example,
in the classic optimal control problem of a falling animal (e.g. a cat) righting
itself in flight, the analysis focuses on the influence of the relative
positions of the front and back halves of the animal on its inertia tensor; the
complex details of the motion are not crucial to a clear understanding of the
problem (see, e.g.,\cite{KS, mont2}, and the references therein). Similar 
situations arise in
computational dynamics where the orientation of a body may be accurately
computed even if the body deformations are not. On the other hand, there are
many situations in which the information of interest is invariant under the
group action. Thus, while in some circumstances it is appropriate to drop to
the base manifold of the principal bundle, in others it may be more natural or
more convenient to remain `upstairs'.

A crucial limitation of the classical theory of connections is the requirement
that the manifold in question be a principal bundle and so the group must act
freely. Many important group actions fail to be free, e.g., transitive actions
on homogeneous spaces, the rotation groups acting on $n\geq 3$-dimensional
Euclidean space, affine and projective actions, etc. Systems with
continuous isotropy arise in a wide variety of applications. In geometric
mechanics, the action of a product group on a manifold of diffeomorphisms or
embeddings, with one factor of the product acting by `body'
transformations, the
other by `spatial' transformations, \cite{L92, L93, L95}, is not free. A
sleeping top (one for which the axis of symmetry of the top is aligned with the
axis of gravity) is a familiar instance of this. In this system, the
infinitesimal versions of `spatial' and `body' rotations (resp. spin and
precession) cannot be distinguished. A more substantial example, which reflects
the same underlying geometric structure, is fluid flow in the absence of external
forces. If the reference region for the fluid mass is axially symmetric, then
any axisymmetric deformation has nontrivial isotropy corresponding to
counter--balancing spatial and body rotations. Such states are analogous to the
sleeping states of the Lagrange top; this common symmetry feature can be
used to
analyze the stability and bifurcation of both systems. 

While one can no longer apply the classical theory of connections on principal
bundles to non-free actions, one would still like to develop a
comparable theory that will carry the many benefits of the usual theory
over to
this context. To this end, we introduce the concept of a {\em
partial connection} that can be applied to general group actions. The key step
is to focus our attention on the connection one-form that can be used to define
the principal bundle connection. The classical Lie algebra-valued connection
form defines an equivariant map from the tangent bundle of the manifold to the
Lie algebra of the transformation group. In our approach, the connection
form is
more appropriately viewed as a projection of each tangent space onto the
infinitesimal group orbit. Such projections can be defined even when the action
fails to be free and lead to promising generalizations of several key
constructs. The shift in focus from generators to projections allows us to
broaden our search for appropriate forms, and we discover that the projections
can be defined using smooth $\fg^*$--valued forms 
even at points in which the isotropy changes.
This projection--based approach is inspired by the reduced energy--momentum
method for stability and bifurcation analysis of relative equilibria of simple
mechanical systems 
\cite{LMSP, SLM} and, in particular, its generalization to non-free actions
 and
regular Lagrangian systems \cite{L92, L93, L95}. In that setting, however, the
decomposition of tangent vectors into their `rigid' and `internal' components
is carried out only at relative equilibria.
Our approach is modeled on the implementation of the simple mechanical
connection by means of the momentum map, and is elucidated in Section 2.

The curvature of a connection can be interpreted as a measure of its
non-integrability, and plays a crucial role in geometric mechanics, mathematical
physics, and differential geometry. In optimal control theory, conditions for
controllability can be described in terms of the curvature of a connection via
the Ambrose-Singer Theorem, 
\cite{AS, ozeki56}. Connection forms facilitate the study of holonomy, which is
important in both classical and quantum mechanics, \cite{ozeki56, berry, SW, 
mont1}. In Section 3, we define a corresponding concept of curvature for a 
partial connections associated with sufficiently regular class of non-free 
actions. In Section 4, we
generalize the classical involutivity result that the curvature is locally zero
if and only if the connection is locally tangent to a cross section of the
principal bundle, and show that a partial connection whose curvature lies
entirely in the tangent spaces of the orbits of a specified isotropy subgroup
is tangent to a slice of the group action. 

Our results have, in part, been directly inspired by the equivariant
approach to
moving frames developed by Fels and Olver \cite{FOI, FOII}. They define a moving
frame as an equivariant map from a manifold
$M$ to a group $G$ acting (locally) freely on 
$M$. These maps, and their linearizations, can be used in the design of
numerical integration schemes. If $M$ is a principal bundle, then the 
trivialized linearization of a moving frame on $M$ is a connection form.
In prior work, \cite{LN00, LO, LN02}, we developed a generalization, 
called {\em partial moving frames}, which play a similar role for
non-free actions. Partial moving frames are mappings from a manifold
$M$ to a group $G$ acting on $M$ that are equivariant modulo isotropy. Like a
genuine connection form, the trivialized linearization of a partial moving frame
can be used to map vector fields on $M$ into trivialized vector fields on
$G$. These
linearizations behave very much like our algebra--valued generalizations of
connection forms, but they typically fail to be equivariant. By relaxing that
condition to relative equivariance, i.e. equivariance with
respect to a specified choice of representatives of the equivalence classes of
the group modulo isotropy, we obtain a further generalization of connections
that includes our motivating examples, and, in fact, served as our original
definition of a partial connection. We describe some of the key features of
relatively equivariant partial connections and the associated forms in 
Section 5.

This work was originally motivated by the need to develop symmetry-preserving
numerical methods (geometric integrators) for solving ODEs and PDEs
with nontrivial symmetry groups. Such methods rely on numerical
integration schemes for initial value problems on Lie groups, cf.
\cite{LS94, LS96, lewis, MKRK, RKMK} and references therein. However, if the
group action has continuous isotropy, the induced vector field on the
group is not uniquely determined. While the true flows of different choices
of vector field will yield the same solution curves back on the manifold,
numerical approximations of these flows will typically yield different
approximate discrete trajectories. Preliminary results in both toy problems and
more substantial micromagnetic calculations are quite encouraging. Adapting the
powerful machinery of connection forms and moving frames to the non-free context
will be of fundamental importance in the further development of such numerical
algorithms.

\section{Partial connections}

We shall assume throughout that a Lie group $G$ acts continuously on a 
manifold $M$. For convenience, we introduce the following notation.
Let $\Phi_g: \cM \to \cM$ and $\apply{m}: G \to \cM$ denote the maps
\[
\apply{m}(g) := \Phi_g(m) := g \cdot m
\]
and, given $\xi \in \fg = T_e G$, 
let $\xi_\cM$ denote the vector field $\ig{\xi}{m} := d \apply{m} \xi$,
called the {\em infinitesimal generator} associated to $\xi$.

A connection on a principal bundle $P$ is a differential system, i.e. a
distribution, $\Gamma$ satisfying $T P = \tfg \oplus \Gamma$ and
$d \Phi_g \cdot \Gamma_p = \Gamma_{g \cdot p}$
for all $p \in P$ and $g \in G$. Here $\tfg$ denotes the differential system
of the tangent spaces to the group orbits, i.e. $\tfg|_p
= T_p (G \cdot p) = \range{d_e \actM p}$. 
Note that we will not explicitly indicate the basepoint of a differential 
structure when it is clear from the context. Specification of a 
connection $\Gamma$ is equivalent to specification of an equivariant 
$\fg$--valued one--form $\alpha$, called the {\em connection form}, satisfying
\[
\alpha \circ d_e \actM p = \mbox{id}, \qquad 
\mbox{i.e.} \quad \alpha_p(\xi_P(p)) = \xi \qquad \mbox{for all $\xi \in \fg$},
\]
for all $p \in P$, where $\alpha_p = \alpha|T_p P$, the restriction of $\alpha$
to the fiber over $p$. By equivariance we mean that $\alpha \circ d \Phi_g
= \Ad g \circ \alpha$ for all $g \in G$. The connection $\Gamma$ and
connection form $\alpha$ are related by the condition
$\mbox{ker} \, \alpha = \Gamma$, i.e. 
$\mbox{ker} \, \alpha_p = \Gamma|_p$ for all $p \in P$. 

We shall retain most the properties of connections and connection forms 
given above in our proposed extension of connections to general actions;
however, our connections need not be differential systems, since the 
dimensions of the group orbits need not be constant throughout the manifold.
The connection form is not uniquely determined by the connection, since
the isotropy algebra $\fg_m = \kernel{d_e \actM m}$ can be nontrivial. However,
equivariant assignments of complements $\parconn$ to $\tfg$ are in one to one 
correspondence with equivariant projections onto $\tfg$; specifically, the 
kernel of the projection is a complementary differential system. As we now 
show, projections onto $\tfg$ are naturally related to $\fg$--valued one forms. 

\begin{Proposition}
\label{conn_relations1}
Given a $\fg$--valued one--form $\alpha$, define the map
$\prpcf := d \actM{} \circ \alpha : T \cM \to \tfg$, i.e.
$\prpcf|T_m M = d \actM m \circ \alpha_m$. 
$\prpcf$ is an equivariant projection onto $\tfg$ if and only if 
$\alpha(\ig{\eta}{m}) = \eta \mod \fg_m$, i.e. 
\beq{proj_range}
\range{(\Id - \alpha \circ d_e \actM m)} \subseteq \fg_m
\eeq
for all $m \in \cM$, and
$\alpha$ is equivariant modulo isotropy, i.e.
\beq{alt_alpha_eq_cond}
\Phi_g^*\alpha_m = \Ad g \alpha_m 
	\quad \mbox{mod} \quad \fg_{g \cdot m}
\eeq
for all $g \in G$ and $m \in \cM$.

Any $\fg$--valued one--form satisfying \for{proj_range}
is discontinuous at singular points of $M$.
\end{Proposition}

\prf
If $\alpha$ satisfies \for{proj_range}, it follows that 
$T M = \tfg \oplus \kernel{\prpcf}$
and $\prpcf|{\tfg} = \Id$; hence $\prpcf$ is a projection onto 
$\tfg$. 
On the other hand, if $\prpcf$ is an equivariant projection onto $\tfg$,
then $\tfg|_m = \range{d_e \actM m}$ for all $m \in M$ and 
$\Id = \prpcf|{\tfg}$ imply \for{proj_range}.

Fix $m \in M$ and $g \in G$. The identity
$(\Ad g \xi)_\cM(g \cdot m) = d_e (\Phi_g \circ \actM m) \xi$
for all $\xi \in \fg$ implies that any $\fg$--valued one--form $\alpha$ 
satisfies
\beqa
d \actM {g \cdot m} \circ \lp \Phi_g^*\alpha - \Ad g \alpha\rp_m
&=& d \actM {g \cdot m} \circ \alpha \circ d_m \Phi_g 
	- d (\Phi_g \circ \actM m) \circ \alpha_m \\
&=& \bbP_\alpha \circ d_m \Phi_g - d_m \Phi_g \circ \bbP_\alpha. 
\eeqa
Thus $\alpha$ is equivariant modulo isotropy if and only if $\prpcf$ is
equivariant.

Given $m \in M$, define the map 
$\pi_m:= \Id - \alpha \circ d_e \actM m: \fg \to \fg$. Since
$\fg_m = \kernel{d_e \actM m}$, we have $\pi_m|\fg_m = \Id$; hence if 
$\alpha$ satisfies \for{proj_range}, then $\range{\pi_m} = \fg_m$.
Continuity of the action implies that $\dim \fg_m \leq \dim \fg_{m_0}$ on a 
sufficiently small neighborhood of a point $m_0$, with equality if and only 
if $m_0$ is regular, and the map $m \mapsto \pi_m$ is continuous if 
$\alpha$ is continuous. Hence if $\alpha$ is continuous at a point $m_0$, then 
\[
\dim \fg_m = \rank{\pi_m} \geq \rank{\pi_{m_0}} = \dim \fg_{m_0}
\]
on a sufficiently small neighborhood of $m_0$. It follows that $\alpha$ can
be continuous at $m_0$ only if $m_0$ is regular.
\prfend

Proposition \ref{conn_relations1} implies that
a $\fg$--valued form $\alpha$ determining a family of equivariant projections 
will be singular at points at which there is a jump in isotropy. Thus, unless 
all points of $M$ are regular, and hence all group orbits in connected
components of $M$ are of equal dimension, 
complements to the tangent spaces to the group orbits
do not appear to have `natural' characterizations in terms of smooth 
$\fg$--valued one--forms. Rather than attempt to specify directly just what
kinds of lapses in smoothness are permissible, we shall depart from the
traditional approach and define our extension of connections utilizing
smooth $\fg^*$--valued forms modeled on momentum maps.
The following example typifies the situations that we will address.

Let $M$ be a Riemannian manifold and let $G$ be a subgroup of the group of 
(local) isometries of $M$. The orthogonal complement 
$\tfg^\perp$ to $\tfg$ in $TM$ will be the prototypical example of
a partial connection. If the action of $G$ is free and proper, and hence $M$
is a principal bundle, the associated connection form is known as 
the {\em simple mechanical connection form} and is given as follows:
Define the equivariant $\fg^*$--valued one--form $\mu$ by
\[
\mu(v) \cdot \xi := \la v, \xi_M(m) \ra_m
\]
for $v \in T_m M$. Note that $\tfg^\perp = \kernel{\mu}$.
The form $\mu$ is the momentum map associated to the $G$ action and the 
Lagrangian $L(v) = \half \norm{v}^2$.
The {\em locked inertia tensor} $\chi: M \to \cL(\fg, \fg^*)$ is given by
\[
(\chi(m) \xi) \cdot \eta = \la \xi_\cM(m), \eta_\cM(m) \ra_m.
\]
The locked inertia tensor is singular precisely at points with
continuous isotropy. Thus if $M$ is a principal bundle, then $\chi(m)$ is 
invertible for all $m$ and $\alpha_m = \chi(m)^{-1} \mu_m$
is the simple mechanical connection form. Note that the $\fg$--valued one--form
$\alpha$ is defined via the $\fg^*$--valued one--form $\mu$ and associated 
map $\chi$; thus, in this setting, the $\fg^*$--valued form can be regarded
as more `directly' linked to the group action and the geometry of the 
problem than the $\fg$--valued form. 
At a point $m$ with continuous isotropy, the relation
$\chi(m) \alpha_m = \mu_m$ does not uniquely determine $\alpha_m$. However,
any $\fg$--valued one--form satisfying $\chi \, \alpha = \mu$ will satisfy
$\tfg^\perp|_m = \alpha_m^{-1}(\fg_m)$ and hence
will qualify as a generalized connection form. Note, however, that in many
situations it appears to be more convenient to work directly with the maps
$\mu$ and $\chi$, since these maps are both smooth, while $\alpha$ will have
singularities at points where the isotropy jumps. 

\exam{$SO(3)$ acting on $\R^3$}
For the sake of concreteness, we now specialize the situation described above 
to the case $M = \R^3$ and $G = SO(3)$, the group of rotations of $\R^3$. 
This action is not free: any nonzero vector $m$ is fixed 
by the circle of rotations about $m$, while the origin is fixed by the entire
group $SO(3)$. The group orbit $G \cdot m$ through $m$ is the sphere of 
radius $||m||$ centered at the origin; 
$\tfg|_m = \lcb \xi \times m : \xi \in \R^3 \rcb$
is the space of infinitesimal rotations of $m$. 
The orthogonal complements to $\tfg$ satisfy
\[
\tfg^\perp|_m = \lcb \begin{array}{ll}
\mbox{span}[m] \quad & m \neq 0 \\
\R^3 & m = 0
\end{array} \right . .
\]
The angular momentum can be regarded as a $so(3)^*$--one--form 
$\mu(v) = m \times v$ for $v \in T_m \R^3$, with $\kernel{\mu} = \tfg^\perp$. 
More generally, given any smooth function $q: \R \to \R$ that is strictly
positive on $\R^+$,
\beq{q_so3}
\mu^q(v) := q(||m||^2) m \times v, 
\eeq
for all $v \in T_m M$, also satisfies $\kernel{\mu^q} = \tfg^\perp$. The 
associated inertia 
tensors $\chi^q(m) := \mu^q \circ d_e \actM m \in \R^{3 \times 3}$ satisfy
\[
\chi^q(m) \xi 
= \lcb \begin{array}{ll}
q(||m||^2) ||m||^2 \proj \perp \xi \quad & m \neq 0 \\ 
0 & m = 0
\end{array} \right . 
\]
where $\proj \perp|T_m M$ denotes orthogonal projection onto 
$\mbox{span}[m]^\perp$. 
Any smooth function $f: \R^3\backslash \{0\} \to \R$ 
determines a $\R^3$--valued one--form 
\beq{so3r3_alpha}
\alpha(v) = \lcb \begin{array}{ll}
||m||^{-2} m \times v + f(m) \la m, v \ra m \quad 
& v \in T_m \R^3, \ m \neq 0 \\ 
0 & v \in T_0 \R^3 
\end{array} \right . 
\eeq
that in turn determines an equivariant projection $\proj \alpha$, with
$\kernel{\proj \alpha} = \tfg^\perp$.
Note that $\alpha$ is discontinuous at the origin for any $f$, while the 
generalized angular momentum $\mu^q = \chi^q \alpha$ is everywhere smooth.
\rmkend

Motivated by the role of the momentum map in our prototypical example, we 
introduce an alternative to $\fg$--valued connection forms; we specify a
differential system complementary to the tangent spaces of the group orbits by
means of a $\fg^*$--valued form that is required to remain 
smooth even at points at which a jump in isotropy occurs. We first show that
an equivariant and appropriately nondegenerate $\fg^*$--valued form determines
a projection onto the tangent spaces of the group orbits, just as the 
connection form does.

\begin{Proposition}
\label{conn_relations3}
An equivariant $\fg^*$--valued one--form $\mu$ satisfies 
$T M = \tfg \oplus \kernel{\mu}$
if and only if the associated equivariant map 
$\chi: \cM \to \cL(\fg, \fg^*)$ given by 
\beq{chi_mu}
\chi(m) := \mu \circ d_e \actM m 
\eeq
satisfies 
\beq{chi_mu_cond}
\kernel{\chi(m)} = \fg_m 
\sands
\range{\chi(m)} = \range{\mu_m}
\eeq
for all $m \in \cM$.

For each $m \in M$, we can 
define the isomorphism $\gamma(m): \range{\chi(m)} \to \tfg$ by 
\beq{gamma_def}
\gamma(m)(\nu) = \xi_M(m) 
\qquad \qquad \mbox{for any $\xi \in \fg$ satisfying $\chi(m) \xi = \nu$}
\eeq
and equivariant projection $\proj \mu := \gamma \circ \mu: TM \to \tfg$. 
\end{Proposition}
\prf
Equivariance of $\mu$ implies that
$\kernel{\mu_{g \cdot m}} = d \Phi_g(\kernel{\mu_m})$ and
\[
\chi(g \cdot m) = \Ads {g^{-1}} \circ \chi(m) \circ \Ad {g^{-1}}
\]
for all $m \in M$ and $g \in G$.
If $\mu$ satisfies $TM = \tfg \oplus \kernel{\mu}$, then 
\[
\kernel{\chi(m)} = \kernel{\mu \circ d_e \actM m} 
= \kernel{d_e \actM m} = \fg_m
\]
and
$\range{\mu_m} = \range{(\mu|{\tfg})|_m} = \range{\chi(m)}$ for all $m \in M$.
 
On the other hand, if $\chi$ given by \for{chi_mu} satisfies \for{chi_mu_cond},
then the nondegeneracy condition $\kernel{\chi(m)} = \fg_m$ implies that
$\mu|{\tfg}$ is injective, while $\range{\chi(m)} =
\range{\mu_m}$ implies that for any $m \in M$ and $v \in T_m \cM$, there exists
$\xi \in \fg$ such that
\[
0 = \mu(v) - \chi(m) \xi = \mu(v - \xi_\cM(m)).
\]
Hence $T \cM = \tfg \oplus \kernel{\mu}$.

The map $\gamma$ is well--defined, since 
$\kernel{\chi(m)} = \kernel{d_e \actM m}$, and equivariant. Equivariance of
$\proj \mu$ thus follows from the equivariance of $\mu$. For any $m \in M$,
\[
\proj \mu \circ d_e \actM m
= \gamma(m) \circ \mu \circ d_e \actM m 
= \gamma(m) \circ \chi(m) 
= d_e \actM m,
\]
so $\proj \mu|\tfg = \Id$. Since $\proj \mu|\kernel{\mu} = 0$ 
and $T \cM = \tfg \oplus \kernel{\mu}$, $\proj \mu$ is a 
projection onto $\tfg$.
\prfend

We now have at hand the appropriate notion of smoothness to formulate our
extension of connections to general actions.

\begin{Definition}
An {\em equivariant partial connection} is a (singular) equivariant
differential system $\parconn$ satisfying 
$T M = \tfg \oplus \parconn$ and either $\parconn$ is smooth or
$\parconn = \kernel{\mu}$ for some smooth $\fg^*$--valued one--form.

An {\em equivariant dual connection form} is a 
smooth equivariant $\fg^*$--valued one--form $\mu$ on $\cM$ satisfying
$T_m M = \tfg \oplus \kernel{\mu_m}$ for all $m \in M$.

An {\em equivariant partial connection form} is 
a $\fg$--valued one--form $\alpha$ on $\cM$ such that the map 
$\prpcf = d \widehat \Phi \circ \alpha$ is an equivariant
projection onto $\tfg$ and $\parconn := \kernel{\prpcf}$
determines an equivariant partial connection. 

An {\em equivariant inertia factor} is an equivariant map 
$\chi: \cM \to \cL(\fg, \fg^*)$ satisfying 
$\kernel{\chi(m)} = \tfg|_m$ for all $m \in M$.
\end{Definition}

From now on, with the exception of \S \ref{beta_equiv}, we shall drop the 
adjective `equivariant' 
and simply refer to partial connections, dual connection forms, etc. 
Note that if $M$ is finite dimensional, $m$ is a regular point of $M$, 
and $\Gamma$ is a partial connection, then there is a neighborhood $\cU$ 
of $m$ such that $\Gamma|\cU$ is a smooth differential system. (Either
$\Gamma$ is a priori smooth or there is a smooth form $\mu$ such that
$\Gamma = \kernel{\mu}$; the rank of $\mu$, and hence the dimension of 
$\Gamma$, is constant on a sufficiently small neighborhood of a regular
point.)

We shall rely primarily on dual connection forms in our calculations. The
following theorem, which is proved in the appendix,
presents some of the fundamental links between dual
connection forms, partial connections, and partial connection forms.

\begin{Theorem}
\label{conn_relations4}
\begin{enumerate}
\item
A dual connection form $\mu$ determines a partial connection 
$\parconn = \kernel{\mu}$ and inertia factor 
$\chi = \mu \circ d_e \actM {}$.

\item
An equivariant (singular) differential system $\parconn$ satisfying
$T M = \tfg \oplus \parconn$ is a partial 
connection if there is inertia factor $\chi$ such that the equivariant
$\fg^*$--valued one--form $\mu$ given by
\beq{pc_variant}
\mu|{\parconn} := 0 
\sands
\mu \circ d_e \actM m := \chi(m)
\quad \mbox{for all $m \in M$} 
\eeq
is smooth, and hence a dual connection form.

\item
A $\fg$--valued one--form $\alpha$ is a partial connection form if 
there is an inertia factor $\chi$ such that $\mu = \chi \, \alpha$ is a dual
connection form with inertia factor $\chi$.
\end{enumerate}
\end{Theorem}

We shall call a pair $(\alpha, \chi)$ consisting of a partial connection form
$\alpha$ and a inertia factor $\chi$ such that 
$\mu := \chi \, \alpha$ is a dual connection form a 
{\em partial connection pair}. Note that, given a partial connection pair,
the inertia factor is not needed to specify the associated partial 
connection; however, it converts a typically singular form to a smooth one, 
and thus will play a crucial role in our development of the curvature form.

\rmk
Note that in contrast to the situation for classical connection forms, 
a partial connection form need not satisfy 
$\kernel{\proj \alpha|T_m M} = \kernel{\alpha_m}$. 
For example, the one--form \for{so3r3_alpha} on $\R^3$ satisfies
\beq{clean}
\range{\alpha_m} \cap \fg_m = \lcb 0 \rcb
\eeq
if and only if $f \equiv 0$. However, given an arbitrary \pcf $\alpha$, we can 
construct a \pcf satisfying \for{clean} and determining the same partial 
connection as $\alpha$. Specifically, set
$\tilde \alpha := \alpha \circ \proj \alpha$. By construction,
$\proj {\tilde \alpha} = \prpcf$
for all $m \in M$. Equivariance modulo isotropy of $\tilde \alpha$ follows 
from that of $\alpha$ and the associated projection; specifically,
$d \Phi_g \circ \proj \alpha = \proj \alpha \circ d \Phi_g$ 
implies that
\[
\range{(\Phi_g^*\tilde \alpha - \Ad g \tilde \alpha)_m}
= \range{(\Phi_g^* \alpha - \Ad g \alpha) \circ \proj \alpha|T_m M}
\subseteq \fg_{g \cdot m}
\]
for any $m \in M$ and $g \in G$. 
\rmkend

\section{Curvature}

The curvature $\Omega$ of a classical connection on a principal bundle
is typically defined as the covariant derivative of the 
unique associated connection form $\alpha$, i.e. 
$\Omega = \cod \alpha = {\mathbb P}^* \rmd \alpha$,
where $\proj {}$ denotes projection onto the connection.
When extending this notion to partial connections, we encounter two 
key difficulties: first, the nonuniqueness of the partial connection form; 
second, the singularity of partial
connection forms at points where there is a jump in isotropy. Ambiguities 
regarding elements in the isotropy subalgebras motivate us to adopt the 
convention that the curvature takes values in the vertical structure $\tfg$, 
rather than in the algebra $\fg$. Complications arising from singularities
are less easily resolved, but we identify a large family of dual connection
forms for which a natural version of curvature can be defined at all points.

Following the standard definition, we define the {\em covariant derivative}
$\cod$ of a (vector--valued) $k$--form $\nu$ with respect to a 
partial connection $\parconn$ as
\[
\cod \nu(v_0, \ldots, v_k) = \rmd \nu(\prpc v_0, \ldots, \prpc v_k)
\]
for any $v_j \in T_m \cM$, $j = 0, \ldots, k$, and $m \in \cM$.
We first consider the situation closest
to the standard one, namely a partial connection form at a regular point.
A partial connection form is smooth at regular points; hence 
we define the {\em curvature} of a partial connection form $\alpha$ 
at a regular point $m$ as
\beq{proto_curv}
\Omega(u, v) = d \actM m(\cod \alpha(u, v))
\eeq
for any vectors $u$, $v \in T_m M$.

To extend our definition of curvature to singular points, we will replace 
partial connection forms, which are discontinuous at singular points. with dual 
connection forms, which are everywhere smooth. If $(\alpha, \chi)$ is a partial 
connection pair associated to a partial connection $\parconn$, we have
\beq{curv_rel}
\chi \, \cod \alpha = \cod (\chi \, \alpha) 
	- \cod \chi \wedge \alpha \circ \proj \parconn
\eeq
at regular points.
The right hand side is well-defined even at points at which $\alpha$ fails
to be differentiable. This suggests that we use \for{curv_rel} to
define the curvature, replacing the dual connection form $\chi \, \alpha$
with a general dual connection form $\mu$ as appropriate.
However, it need not be the case that the right hand side of \for{curv_rel}
lie in the range of $\chi(m)$. Before addressing this issue, we show that
the wedge product appearing in \for{curv_rel} is identically zero, yielding
$\chi \, \cod \alpha = \cod (\chi \, \alpha)$ whenever the left hand side is 
defined.

\begin{Lemma}
\label{good_chi}
Let $\mu$ be a dual connection form with inertia factor $\chi$.
Then $\range{\cod \chi(m)} \subseteq \Ann \fg_m$, the annihilator of $\fg_m$. 
\end{Lemma}
Proof: 
Equivariance of $\mu$ implies that $\eta_M(\mu) = - \ads \eta \mu$
for any $m \in M$ and $\eta \in \fg$. 
Hence $\mu(\eta_M) = \chi \, \eta$ implies that
\beq{diff_eq_1}
\interior {\eta_M} \rmd \mu
= \eta_M(\mu) - \rmd (\mu(\eta_M))
= - \ads \eta \mu - \rmd \chi \, \eta. 
\eeq
If $\zeta \in \fg_m$, and hence $\zeta_M(m) = 0$, 
and $u \in \parconn|_m = \kernel{\mu_m}$, then \for{diff_eq_1} yields
$0 = \rmd \chi(u) \zeta$.
\prfend

Combining \for{proto_curv}, \for{curv_rel}, and Lemma \ref{good_chi}, given a
dual connection pair $(\alpha, \chi)$, we see that at a regular point $m$
\[
\Omega_m = d \actM m(\codm \alpha) 
= \gamma(m)(\chi(m) \codm \alpha) 
= \gamma(m)(\codm (\chi \, \alpha)).
\]
Note that the right hand side is well--defined even at singular points, provided
that 
\[
\range{\codm (\chi \, \alpha)} \subseteq \range{\chi(m)}.
\]
Motivated by this observation, we introduce the following definitions:

\begin{Definition}
A dual connection form $\mu$ is {\em docile} at $m$ if
$\range{\codm \mu} \subseteq \range{\mu_m}$.
If $\mu$ is docile at $m$, then the {\em curvature} of $\mu$ at $m$ is 
$\Omega_m := \gamma(m) \circ \cod \mu_m$,
where $\gamma$ is given by \for{gamma_def}.

If a partial connection form $\alpha$ is differentiable at $m$, then the 
{\em curvature} of $\alpha$ at $m$ is given by \for{proto_curv}.
A partial connection pair $(\alpha, \chi)$ is {\em docile} 
at $m$ if the dual connection form $\chi \, \alpha$ is docile at $m$.
If $(\alpha, \chi)$ is docile at $m$, then the {\em curvature} 
of $(\alpha, \chi)$ at $m$ equals that of $\chi \, \alpha$.
\end{Definition}

In the classical setting, there is a unique connection form associated to a
given connection. Thus the curvature can naturally be regarded as data of
the connection itself, not just of the connection form. 
In general, there is not a unique partial connection form or dual connection 
form associated to a given partial connection; hence curvature is not a priori 
determined by the partial connection, rather than a specific form. 
We shall show in \S 5 that at regular points of $M$, i.e. points at
which the dimensions of the group orbits are locally constant, the curvature 
of a partial connection is well--defined. However, at
singular points it can easily occur that two dual connection forms 
determining the same partial connection fail to share docility. 
For example, the dual connection forms \for{q_so3} for the action of 
$SO(3)$ on $\R^3$ satisfy
\[
{\rm d}\mu^q(0)(u_0, v_0) = 2 \, q(0) \, u_0 \times v_0,
\]
while $\range{\mu^q(0)} = \lcb 0 \rcb$. If $q(0) \neq 0$, then 
${\rm d}\mu^q(0)(u_0, v_0) \not \in \range{\mu^q(0)}$ for any $u_0$ and $v_0$ 
that are not parallel.
Thus the dual connection forms $\mu^q$ discussed above are docile if 
and only if $q(0) = 0$, in which case the curvature at the origin is zero.
This example suggests that in some circumstances we can `tame' a given dual 
connection form that fails to be docile, obtaining a docile dual connection 
form whose curvature agrees with that of the original form wherever the 
original form is docile. 
 
\begin{Remark}
\label{tame_conn}
{\rm
Consider two inertia factors $\chi$ and $\tilde \chi$ compatible with a partial
connection $\Gamma$, in the sense that \for{pc_variant} determines smooth 
forms $\mu$ and $\tilde \mu$, and satisfying
$\tilde \chi = \sigma \circ \chi$ for some smooth map 
$\sigma: M \to \cL(\fg^*, \fg^*)$. 
If $\Omega$ (respectively $\tilde \Omega$) denotes the curvature 
of $\mu$ (respectively $\tilde \mu$), then 
\[
\cod \tilde \mu = \sigma \, \cod \mu 
	+ \cod \sigma \wedge (\mu \circ \proj \parconn {}{})
	= \sigma \, \cod \mu 
\]
and $\tilde \gamma \circ \sigma = \gamma$ imply that
\[
\tilde \Omega
= \tilde \gamma \circ \cod \tilde \mu 
= \gamma \circ \cod \mu 
= \Omega
\]
wherever both $\mu$ and $\tilde \mu$ are docile. 

If $\tilde \mu$'s domain of docility is larger than that of $\mu$, it
may be preferable to replace $\mu$ by $\tilde \mu$.
For example, assume that $\fg$ is a inner product space and $\mu$ is a dual 
connection form such that $\chi(m)$ is symmetric for every $m \in M$. 
Let $\mu^\sharp$ denote the $\fg$--valued one--form satisfying
$\la \mu^\sharp(v), \xi \ra = \mu(v) \cdot \xi$ for all $\xi \in \fg$ and 
$v \in TM$; the dual connection form $\tilde \mu = \chi \circ \mu^\sharp$ is 
docile on all of $M$. Since $\mu$ and $\tilde \mu$ determine the same 
partial connection and have the same curvature wherever the curvature 
of $\mu$ is defined, the form $\tilde \mu$ is, by some standards, the 
preferable one to use when analysing the partial connection. This situation 
frequently arises when $\mu$ is the momentum map determined by the kinetic 
energy on a Riemannian manifold $\chi$ is the `locked inertia tensor', as 
in our prototypical example. We investigate such an example in the next
section.
}
\prfend
\end{Remark}

Using the equivariance properties of dual connection forms, we can derive
the analogs of the classical structure equations for dual connection forms
and partial connection pairs.

\begin{Proposition}
\label{structure}
Given a dual connection form $\mu$ that is docile at $m$ and tangent vectors 
$u$, $v \in T_m M$, let $\xi$, $\eta \in \fg$ satisfy 
$\mu(u) = \chi(m) \xi$ and $\mu(v) = \chi(m) \eta$. Then
\[
\Omega(u, v) + [\xi, \eta]_M(m) = \gamma(m)(\rmd \mu(u, v) 
	- \rmd \chi(u) \eta + \rmd \chi(v) \xi).
\]

Given a partial connection pair $(\alpha, \chi)$, define the $\fg$--valued 
curvature form $\Omega^\alpha := \alpha \circ \Omega$ at points where 
$(\alpha, \chi)$ is docile. Then
\[
\chi (\Omega^\alpha + \alpha \wedge \alpha) = \rmd (\chi \, \alpha) 
	- \rmd \chi \wedge \alpha,
\]
where $(\alpha \wedge \alpha)(u, v) = [\alpha(u), \alpha(v)]$
for all $u$, $v \in TM$. In particular, if $\rmd \chi(m) = 0$, we 
recapture the classical structure equations modulo $\fg_m$ at $m$.
\end{Proposition}
\prf
Equivariance of $\chi$ implies that
$\xi_M(\chi) = - \ads \xi \circ \chi - \chi  \circ \ad \xi$
for all $\xi \in \fg$. Thus \for{diff_eq_1} implies that 
\[
\rmd \mu(\xi_M, \eta_M)
= - \ads \xi (\chi \eta) + \ads \eta (\chi \xi) - \chi[\xi, \eta].
\]
and hence
\beqa
\cod \mu(u, v) &=& \rmd \mu(u - \xi_M(m), v - \eta_M(m)) \\
&=& \rmd \mu(u, v) - \rmd \chi(u) \, \eta 
	 + \rmd \chi(v) \, \xi - \chi(m)[\xi, \eta]. \qquad \qquad \prfend
\eeqa

\subsection{An example: the combined left--right action of a subgroup on a
Lie group}
\label{subgrp_example}
\newcommand{\Adgi}{\mbox{Ad}_{g^{-1}}}
\newcommand{\brackxo}{[\xi', \omega']}

Let $G$ be a Lie group, with Lie algebra $\fg$ equipped with an 
$\Ad{}$--invariant inner product $\la \ , \ \ra$, and let $H$ be a subgroup 
of $G$. Let $H \times H$ act on $G$ by 
\[
(h, k) \cdot g = h g k^{-1}.
\]
Given a subspace $V \subseteq \fg$, let $V^\perp$ denote its orthogonal 
complement and $\proj V {}: \fg \to V$ the orthogonal projection onto $V$.
We take as our partial connection $\Gamma = (\fh \times \fh)^\perp$, the 
orthogonal complement to the tangent to the group orbit with respect to the 
right--invariant metric induced on $G$ by $\la \ , \ \ra$. If we identify 
$T_g G$ with $\fg$ by right trivialization, as we shall throughout this 
example, then the infinitesimal generator of 
$(\eta, \zeta) \in \fh \times \fh$ is 
\[
(\eta, \zeta)_G(g) = \eta - \Ad g \zeta.
\]
If we let $\Ad g(\fh)$ denote the image of $\fh$ under the adjoint action 
of $g$, it follows that the isotropy subalgebra of $g$ is  
\[
(\fh \times \fh)_g = (\Id \times \Adgi)(\fh \cap \Ad g(\fh)).
\]

The momentum map $\mu$ associated to the action of $H \times H$ on $G$ and
the Lagrangian $L(\xi) = \half |\xi|^2$ is a dual connection form compatible
with $\Gamma = \widetilde {(\fh \times \fh)}^\perp$. However, as we shall 
show, $\mu$ fails to be docile at points in the normalizer $N(H)$ of $H$. Hence 
we work with the `tamed' dual connection form $\nu$ associated to $\mu$, as
described in Remark \ref{tame_conn}. The form $\nu$ has curvature
\beq{curv_ex}
\Omega_g(\xi, \omega) 
= (\projeh - \projegh)[\proj{\,\Gamma_{\! g}} \xi, \proj {\, \Gamma_{\! g}} \eta]
\eeq
for all $g \in G$ and $\xi$, $\omega \in \fg$. It follows that
\[
\mbox{rank}\, \Omega_g 
\leq 2 (\mbox{dim}\, \fh - \mbox{dim} (\fh \cap \Ad g(\fh))).
\]
In particular, $\Omega_g = 0$ if $g \in N(H)$.

If we identify $\fh^*$ with $\fh$ using $\la \ , \ \ra$, then the momentum
map $\mu$ satisfies
\[
\mu_g(\xi) = (\projeh \xi, - \projeh \Adgi \xi)
= (\projeh \xi, - \Adgi \projegh \xi),
\]
with associated correction factor
\[
\chi(g) = \lp \begin{array}{cc} \Id & - \projeh \Ad g \\
-  \projeh \Ad {g^{-1}} & \Id \end{array} \rp.
\]
We first compute the exterior derivative of $\mu$. If we let $X_\xi$ denote the 
right invariant vector field associated to $\xi \in \fg$, then 
\[
X_\xi(\mu(X_\omega))(g) = (0, \projeh \Adgi [\xi, \omega])
\]
and $[X_\xi, X_\omega] = - X_{[\xi, \omega]}$ imply that
\[
\rmd \mu_g(\xi, \omega) 
= (X_\xi(\mu(X_\omega)) - X_\omega(\mu(X_\xi)) - \mu([X_\xi, X_\omega]))(g) 
= (\projeh {} [\xi, \omega], \projeh \Adgi [\xi, \omega]).
\]
Note that $\mbox{range}[\cod \mu_g] \subseteq \mbox{range}[\chi(g)]$ 
need not hold for all $g \in G$. In particular, if $g \in N(H)$, then 
$\projeh$ commutes with $\Adgi$; hence, in this case,
\[
\rmd \mu_g(\xi, \omega) 
= (\Id \times \Adgi) \projeh [\xi, \omega]
\]
for all $\xi$, $\omega \in \fg$, while $\range{\chi(g) 
= (\Id \times (-\Adgi)) \fh}$. Thus $\mu$ fails to be docile at points in 
$N(H)$ unless $\fh^\perp$ is a subalgebra of $\fg$. Hence we replace $\mu$ with 
the `tamed' dual connection form $\nu = \chi \circ \mu$, as in Remark
\ref{tame_conn}. Given our
identification of $\fh$ and $\fh^*$, the map $\sharp$ is trivial.

We now show that the curvature $\Omega$ of $\nu$ satisfies \for{curv_ex}.
Given $\xi$, $\omega \in \fg$, let 
\[
\brackxo = \alpha + \beta + \gamma + \delta,
\]
where $\alpha \in \fh \cap (g \cdot \fh)^\perp$, 
$\beta \in \fh^\perp \cap g \cdot \fh$, \quad 
$\gamma \in \fh \cap g \cdot \fh$, and $\delta \in \Gamma_g$, and let
$\xi'$ and $\omega'$ denote the projections of $\xi$ and $\omega$ into
$\Gamma|_g$. Then
\[
(\projeh - \projegh)\brackxo
= (\alpha + \gamma) - (\beta + \gamma) 
= \alpha - \beta.
\]
Observe that
\[
\mu_g(\alpha - \beta) = (\projeh (\alpha - \beta),  
	- \Adgi \projegh (\alpha - \beta))
= (\alpha, \Adgi \beta),
\]
while
\[
\cod \mu_g(\xi, \omega) 
= (\projeh \times \Adgi \projegh) \brackxo
= (\alpha + \gamma, \Adgi (\beta + \gamma)).
\]
Thus
\[
\cod \mu_g(\xi, \omega) - \mu_g(\alpha - \beta) 
= (\gamma, \Adgi \gamma)
\in \kernel{\chi(g)},
\]
and hence
\[
\cod \nu_g(\xi, \omega) 
= \chi(g) \cod \mu_g(\xi, \omega) 
= \chi(g) \mu_g(\alpha - \beta) 
= \nu_g(\alpha - \beta).
\]

As a simple application of the formulas derived above, we now provide an
example of a point with nontrivial isotropy and nonzero curvature. We take
as our manifold the group $SU(3)$ and select as our subgroup $H$ the two--torus 
of diagonal matrices in $SU(3)$. We work with the orthogonal basis 
$\{\delta_1, \delta_2, \sigma_1, \sigma_2, \sigma_3, \xi_1, \xi_2, \xi_3\}$,
where $\delta_1 = \mbox{diag}(i, -i, 0)$, 
$\delta_2 = \mbox{diag}(i, i, -2i)$, $\sigma_j$ has $i$ 
in the $k \ell$ and $\ell k$ positions and zeroes elsewhere, and
$\xi_j$ has 1 in the $k \ell$ position, $-1$ in the $\ell k$ position,
and zeroes elsewhere; here $j k \ell$ is a cyclic permutation 
of 123. 
Note that $\{\delta_1, \delta_2\}$ is a basis for the two--torus 
$H$, while $\{\xi_1, \xi_2, \xi_3\}$ is a basis for the rotation group $SO(3)$. 

We compute the curvature at an element $g$ of $SO(3)$ corresponding 
to a rotation through an angle $\theta$ about the vertical axis, 
$\theta \neq \frac {n \, \pi} 2$ for any $n \in {\mathbb Z}$. The 
adjoint action of $g$ fixes $\delta_2$ and maps $\delta_1$ to 
$\cos 2 \theta \, \delta_1 + \sin 2 \theta \, \sigma_3$. Hence
\[
(\fh \times \fh)_g = \spn{\delta_2}, \quad
\widetilde{(\fh \times \fh)}|_g = \spn{\delta_1, \delta_2, \sigma_3},
\quad \mbox{and} \quad
\Gamma|_g = \spn{\sigma_1, \sigma_2, \xi_1, \xi_2, \xi_3}.
\]
The commutators of the basis elements of $\Gamma|_g$ are
$[\sigma_1, \sigma_2] = \xi_3$, 
$[\xi_j, \xi_k] = \xi_\ell$, where $j k \ell$ is a 
cyclic permutation of 123, and
\[
[\sigma_j, \xi_k] = \lcb \begin{array}{ll}
(-1)^j \delta_1 + \delta_2 \qquad & j = k \neq 3 \\
- \sigma_3 & 3 \neq j \neq k \neq 3 \\
(-1)^j \sigma_{j'} & j \neq k = 3
\end{array} \right .,
\]
where $jj' = 12$ or 21. 
We have $(\proj{\fh} - \proj{\Ad g \fh}) \eta 
= \half (\la \eta, \sigma_3 \ra \sigma_3 - \la \eta, \delta_1 \ra \delta_1)$
for all $\eta \in su(3)$. Hence the nontrivial elements of the curvature 
at $g$ are
\[
\Omega_g(\sigma_j, \xi_k) = \lcb \begin{array}{ll}
(-1)^j \delta_1 \qquad & j = k \neq 3 \\
- \sigma_3 & 3 \neq j \neq k \neq 3 
\end{array} \right .
\]
and the curvature at $g$ has rank two, with
$\range{\Omega_g} = \spn{\delta_1, \sigma_3}$. 

\section{Curvature and involutivity}
\label{almost_hor}

The classical
formula $\Omega(X, Y) = \alpha([Y, X])$ for
the curvature $\Omega$ of a connection form $\alpha$, where $X$ and $Y$ are
horizontal vector fields, can be rephrased as the assertion
that the curvature of a connection measures the extent to which the connection
is involutive. We can easily show that at regular points of $M$ with respect to
proper actions the curvature of a partial connection does not depend on the
choice of partial or dual connection form used to characterize the differential 
system and this classical involutivity relation is satisfied. The development
of an analogous expression for the curvature at singular points involves the
construction of a smooth differential system, the system of `almost horizontal'
vectors, containing the singular partial connection 
on a neighborhood of a singular point. 

The curvature at a regular point can be expressed in terms of the Lie bracket
of horizontal vector fields, as in the classical case. This expression leads
directly to a correspondence between (locally) zero curvature and involutivity 
of the partial connection $\Gamma$, and hence the existence of maximal integral 
manifolds tangent to $\Gamma$. 

\begin{Theorem}
\label{reg_curv}
Let $\Omega$ denote the curvature of a dual connection form or partial 
connection form compatible with the partial connection $\Gamma$.
Given any horizontal vector fields $X$ and $Y$, the vector field 
$\Omega(X, Y) + [X, Y]$ is horizontal.

If $\Omega$ is identically zero on some open set containing only regular
points, then $\Gamma$ is involutive on that set.
\end{Theorem}
\prf
Let $Z = [X, Y]$.
We first consider a dual connection form $\mu$. We have
\[
\cod \mu(X, Y)
= \rmd \mu(X, Y)
= X (\mu(Y)) - Y (\mu(X)) - \mu([X, Y])
= - \mu(Z)
\]
and hence 
\[
Z + \Omega(X, Y) 
= Z - \gamma(\mu(Z))
= (\Id - \proj \mu) Z
= \prpc Z.
\]
Given a partial connection form $\alpha$, analogous arguments show that
\[
Z + \Omega(X, Y) 
= Z - d \actM m(\alpha(Z))
= (\Id - \proj \alpha) Z
= \prpc Z. \qquad
\]

As was previously discussed, if $m$ is a regular point, then $\parconn$ is a 
smooth differential system near $m$ and hence is spanned by horizontal vector 
fields on some neighborhood of $m$. Thus the curvature on this neighborhood 
is determined by the relation $\Omega(X, Y) = (\prpc - \Id)[X, Y]$ for 
horizontal vector fields $X$ and $Y$. In particular, if $\Omega$ is 
identically zero on such a neighborhood, then $\Gamma$ is involutive there.
\prfend

The key feature of the proof of Theorem \ref{reg_curv}
is the existence of vector fields 
generating $\Gamma$ near $m$. If $m$ is a singular point, then the dimension of
$\Gamma|_m$ is greater than that of $\Gamma|_n$ at nearby
points $n$. Thus we cannot find a set of horizontal vector fields on a 
neighborhood of $m$ that span $\Gamma|_m$, and hence cannot directly invoke the 
identity relating exterior derivatives and commutators of vector fields. 
At singular points
we must relax the notion of horizontality, obtaining a condition that yields 
a smooth differential system even in the neighborhood of a point at which a
jump in isotropy occurs. We enlarge the partial connection by including a
`rotated' copy of the tangent space to the $\gmo$ orbit at each point. These
spaces are rotated so as to yield a smooth differential system near the given 
singular point. By analogously rotating the dual connection form, we obtain a 
form whose kernel is the desired differential system. Using this adapted form,
we can mimic the argument of Theorem \ref{reg_curv} at singular points under 
an appropriate hypothesis on the adapting map. 

\begin{Definition}
An {\em adaptor} $\phi$ for $m_0$ is a smooth $\gmo$--equivariant map 
$\phi: \cU \to G$ on a $\gmo$--invariant neighborhood $\cU$ of $m_0$ in $M$ 
satisfying $\Ad {\phi(m)} \fgmo \supseteq \fg_m$ for all $m \in \cU$, and 
$\phi(m_0) \in \gmo$.
Here $\gmo$--equivariance is with respect to conjugation in $G$, so that
$\phi(g \cdot m) = g \, \phi(m) g^{-1}$ for all $m \in S$ and $g \in \gmo$.

Given a singular point $m_0$, a partial connection $\Gamma$, and 
an adaptor $\phi$ for $m_0$, $\Xi := \Gamma \oplus \Ad \phi \fgmo$
is the {\em almost horizontal} differential system for $m_0$.

Fix a dual connection form $\mu$ for the partial connection $\Gamma$ and a 
$\gmo$--equivariant projection $\pi$ on $\fg$ with $\kernel{\pi} = \fgmo$.
Define the $\gmo$--equivariant map $\chi_\phi: \cU \to \cL(\fg, \fg^*)$,
where $\cU$ is the domain of the adaptor $\phi$, by
\[
\chi_\phi(m) := \chi(m) \circ \Ad {\phi(m)}.
\]
Using $\chi_\phi$, define the projections
$\pi_\phi(m) := \chi_\phi(m) \circ \pi \circ \chi_\phi(m)^{-1}$ and the 
{\em adapted dual connection form} $\wtm_m := \pi_\phi(m) \circ \mu_m$.
\end{Definition}

Note that $\chi_\phi(m)^{-1}(\mu(v))$ is defined modulo 
$\Ad {\phi(m)^{-1}}(\fg_m) \subseteq \fgmo = \kernel{\pi}$, 
and hence $\pi_\phi$ and $\wtm$ are well--defined. The adapted dual connection 
form $\wtm$ is a constant rank, $\gmo$--equivariant $\fg^*$--valued one form  
with $\kernel{\wtm} = \Xi$. Hence $\Xi$ is a smooth differential system.
Finally, note that $\pi_\phi \circ \chi_\phi = \chi_\phi \circ \pi$.

\begin{Proposition}
\label{adapt_exist}
If $G$ acts properly, then there is an adaptor for any point in $M$. Given a
partial connection $\Gamma$ and a point $m_0 \in M$, there is an adaptor 
$\phi$ for $m_0$ satisfying $\phi(m_0) = e$ and 
$d \phi(\Gamma|_{m_0}) \subseteq \fgmo$. 
\end{Proposition}

The proof of Proposition \ref{adapt_exist}, which makes use of some technical
results related to slices, is given in the appendix. 

\rmk
We could work with the {\em adapted connection form}
$\widetilde \alpha(v) := \Ad {\phi(m)} \pi \, \chi_\phi(m)^{-1}(\mu(v))$,
rather than the adapted dual connection form $\wtm$.
Either $\wtm$ or $\widetilde \alpha$ can be used to define a
$\gmo$--equivariant projection $\proj \phi$ on $T \cU$, with 
$\kernel{\proj \phi} = \Xi$, by 
$\proj {\phi} := \gamma \circ \wtm$ or $\proj {\phi} 
:= d_e \actM \circ \widetilde \alpha$. 
\rmkend

We now generalize the `lack of involutivity' characterization of curvature 
given in Theorem \ref{reg_curv} to singular points, replacing the partial
connection $\Gamma$ with the almost horizontal differential structure $\Xi$.

\begin{Theorem}
\label{sing_comm}
If $\mu$ is docile and $\phi$ is an adaptor for $m_0$ satisfying
\beq{weird_cond}
[\trivl \phi(\Xi|_m), \fgmo] \subseteq \fgmo
\eeq
where $\trivl_m \phi = d_m (L_{\phi(m)^{-1}} \circ \phi)$ is 
the left trivialization of the linearization of $\phi$ at $m$,
for some $m$ near $m_0$, then for any almost horizontal vector fields $X$ 
and $Y$, $\Omega(X, Y) + [X, Y]$ is almost horizontal at $m$.

If \for{weird_cond} holds for all $m$ on a neighborhood $\cV$ of $m_0$ and 
$\range{\Omega} \subseteq \Xi$ on $\cV$, then $\Xi$ is involutive on $\cV$. 
\end{Theorem}

The key ingredient in the proof of Theorem \ref{sing_comm} is an 
expression for the differential of the adapted form $\wtm$ in 
terms of the covariant derivative of the dual connection form $\mu$.
We assume here, and throughout the remainder of this section, that 
the adapted inertia factor $\chi_\phi$ has a differentiable restricted
pseudo--inverse $\iota$; specifically, that there is a differentiable
map $\iota: \cU \to \cL(\fg^*, \fg)$ such that
\[
\pi = \pi \circ \iota(m) \circ \chi_\phi(m)
\]
for all $m \in \cU$. (Note that $\kernel{\chi_\phi(m)} 
= \Ad {\phi(m)^{-1}} \fg_m \subseteq \fgmo = \kernel{\pi}$.)
If we set 
\[
\wtpp(m) := \chi_\phi(m) \circ \pi \circ \iota(m)
\]
for all $m \in \cU$, then $\pi_\phi(m) = \wtpp(m)|{\range{\chi_\phi(m)}}$, 
and hence $\wtm = \wtpp \circ \mu$.   

\begin{Lemma}
\label{rel_curv}
Let $u, v \in T_m \cU$,
and let $\xi$ and $\eta \in \fg$ satisfy $\chi_\phi(m) \xi = \mu(u)$ and
$\chi_\phi(m) \eta = \mu(v)$.
If $\mu$ is docile at $m$ and $\phi$ satisfies \for{weird_cond}, then
\beqa
d \wtm(u, v) - \pi_\phi(m) \cod \mu(u, v)
&=& d \chi_\phi(u) \pi \eta - d \chi_\phi(v) \pi \xi \\
&& \quad {} + \chi_\phi(m) \pi ([\trivl \phi(v), \xi] 
	- [\trivl \phi(u), \eta] + [\xi, \eta]).
\eeqa
\end{Lemma}
\prf
Leibniz's Rule implies that
\beq{first_exp}
d \wtm(u, v) = d \wtpp(u) \mu(v) - d \wtpp(v) \mu(u) + \wtpp(m) d \mu(u, v).
\eeq
Linearizing $\pi = \pi \circ \iota \circ \chi_\phi$ yields
\[
0 = \pi \circ (d \iota(v) \circ \chi_\phi(m)  
	+ \iota(m) \circ d \chi_\phi(v)).
\]
Thus
\beqa
d \wtpp(v) \mu(u) &=& d \wtpp(v) \chi_\phi(m) \xi \\
&=& (d \chi_\phi(v) \pi \iota(m)
	+ \chi_\phi(m) \pi d \iota(v)) \chi_\phi(m) \xi \\
&=& d \chi_\phi(v) \pi \xi - \wtpp(m) d \chi_\phi(v) \xi.
\eeqa
Next, we have
\[
d \chi_\phi(v) \xi = d \chi(v) \Ad {\phi(m)} \xi 
	+ \chi_\phi(m) [\trivl \phi(v), \xi].
\]
Entirely analogous expressions hold when $u$ and $v$ are exchanged.
Proposition \ref{structure} implies that
\[
d \mu(u, v) 
= \cod \mu(u, v) + d \chi(u) \Ad {\phi(m)} \eta 
	- d \chi(v) \Ad {\phi(m)} \xi + \chi_\phi(m) [\xi, \eta].
\]
Substituting these expressions into \for{first_exp}
and regrouping terms yields the desired expression.
\prfend

\prfof{Theorem \ref{sing_comm}}
Let $X$ and $Y$ be almost horizontal vector fields. We apply 
Lemma \ref{rel_curv}, taking $u = X(m)$ and $v = Y(m)$. The condition
$u, v \in \Xi|_m$ implies that $\xi, \eta \in \fgmo = \kernel{\pi}$, and hence
$\pi_\phi \cod \mu(u, v) = d \wtm(u, v)$.
Thus $\wtm(X) = \wtm(Y) = 0$ implies
\[
0 = d \wtm(X, Y) + \wtm([X, Y])
= \pi_\phi (\cod \mu(X, Y) + \mu([X, Y])).
\]
Since $\kernel{\pi_\phi} = \chi_\phi(\fgmo)$, calculations analogous to
those used in the proof of Theorem \ref{reg_curv} yield
\[
\Omega(X, Y) + [X, Y] = \proj \Gamma [X, Y] \qquad \mod \widetilde{\Ad \phi \fgmo}
= \widetilde \fg \cap \Xi;
\]
hence $\Gamma \subseteq \Xi$ implies that $\Omega(X, Y) + [X, Y]$ takes 
values in $\Xi$.

Involutivity of $\Xi$ when $\range{\Omega} \subseteq \Xi$ follows immediately
from the first part of the theorem. 
\prfend

\begin{Cor}
If $G$ acts properly, the curvature at a singular point $m_0$ is determined 
by the equation $\Omega(X, Y)(m_0) = (\proj \Gamma - \Id) [X, Y](m_0)$ 
for almost horizontal vector fields $X$ and $Y$. 
\end{Cor}
\prf
Proposition \ref{adapt_exist} guarantees the existence of an adaptor $\phi$ 
for $m_0$ satisfying \for{weird_cond} at $m_0$. 
Since $\Xi$ is a smooth differential structure spanned by almost horizontal
vector fields, Theorem \ref{sing_comm} determines the curvature modulo 
$\widetilde{\Ad \phi \fgmo}$ at points near $m_0$ satisfying \for{weird_cond}. 
Hence, since $\widetilde{\Ad \phi \fgmo}|_{m_0}$ is trivial, the curvature at 
$m_0$ is entirely determined by Theorem \ref{sing_comm}.
\prfend

We now show that the classical result that the curvature is identically zero
on some neighborhood of a given point if and only if the horizontal 
differential system is tangent to a local cross section through that point
can be generalized to partial connections under appropriate hypotheses on the 
isotropy subgroups. By a local cross section, we mean a submanifold $S_0$ 
such that $G \cdot S_0$ contains a neighborhood of $m_0$ in $M$ and 
$g \cdot m \in S_0$ for $m \in S_0$ only if $g \cdot m = m$. 

\begin{Cor}
\label{flat_reg_lcs}
If $G$ acts properly on $M$, 
$\Omega$ equals zero near $m_0$, and there is an adaptor $\phi$ satisfying 
\beq{full}
G_m = \phi(m) \gmo \phi(m)^{-1} \qquad \mbox{for all $m$ near $m_0$},
\eeq
then $\Gamma$ is tangent to a local cross section through $m_0$.
\end{Cor}

\noindent
The proof of Corollary \ref{flat_reg_lcs} is rather technical and is in part
modeled on a proof for a similar result for slices given in \cite{DK}. Hence 
it is relegated to the appendix. 

\rmk
If $G_m = \exp(\fg_m)$ for all $m$ in a neighborhood of a regular point $m_0$
of a proper action, then equivariance of the exponential map implies that any 
adaptor $\phi$ for $m_0$ satisfies \for{full}.
\rmkend

A slice generalizes the notion of a local cross section, allowing some overlap 
of the slice and the group orbits near singular points. Specifically, a slice 
at $m_0$ is a submanifold $S$ through $m_0$ satisfying
\begin{enumerate}
\item
$T_{m_0} M = T_{m_0} S \oplus \tfg|_{m_0}$ and 
\item
$T_m M  = T_m S + \tfg|_m$ for all $m \in S$
\item
if $m \in S$ and $g \in G$, then $g \cdot m \in S$ if and only if 
$g \in G_{m_0}$.
\end{enumerate}
\begin{Cor}
\label{Abel_curv}
If $G$ acts properly, $\gmo$ is a normal subgroup of $G$, $G_m \subseteq \gmo$ 
for all $m$ in a neighborhood of $m_0$, and 
$\range{\Omega} \subseteq \widetilde {\fgmo}$
on that neighborhood, then some neighborhood of $m_0$ in the integral
submanifold of $\Xi$ containing $m_0$ is a slice.
\end{Cor}

The following lemma establishes the essential information regarding slices 
used in the proof of Corollary \ref{Abel_curv}. 
These results are minor variations of standard results and for 
the most part we utilize straightforward modifications of the proofs given 
in \cite{DK}; hence the proof of the Lemma is given in the appendix.

\begin{Lemma}
\label{slice_cond}
\begin{enumerate}
\item
Given an involutive $\gmo$--equivariant differential system 
$\Delta \supseteq \widetilde {\fgmo}$ on a $\gmo$--invariant neighborhood 
of $m_0$, there
is an $\gmo$--invariant integral submanifold of $\Delta$ containing $m_0$. 
\item
If $G$ acts properly and $S$ is a $\gmo$--invariant submanifold of $M$ containing 
$m_0$ and satisfying
$T_{m_0} M = T_{m_0} S \oplus \tfg|_{m_0}$, there is a neighborhood   
$S_0$ of $m_0$ in $S$ such that $S_0$ is a slice through $m_0$.
\end{enumerate}
\end{Lemma}

\prfof{Corollary \ref{Abel_curv}}
Since $G$ is normal and $G_m \subseteq \gmo$ for all $m$ near $m_0$,
we can use the trivial adaptor $\phi \equiv e$. Theorem \ref{sing_comm}
implies that $\Xi = \Gamma \oplus \widetilde {\fgmo}$ is involutive near $m_0$.

Let $S'$ denote the maximal integral manifold of the almost horizontal system
containing $m_0$. Lemma \ref{slice_cond}.i and the $\gmo$--equivariance of 
$\Xi$ 
imply that $S'$ is $\gmo$--invariant. Lemma \ref{slice_cond}.ii and
\[
T_{m_0} M = (\Gamma \oplus \tfg)|_{m_0}
= T_{m_0} S' + \tfg|_{m_0}
\]
imply that some $\gmo$--invariant
neighborhood $S$ of $m_0$ in $S'$ is a slice through $m_0$.
\prfend

\noindent
{\bf Example: $S^1 \times S^1$ acting on $SO(3)$} 
\medskip

As an application of Corollary \ref{Abel_curv}, we consider a special case of 
\S \ref{subgrp_example}, with $G = SO(3)$ and $H \approx S^1$ consisting of rotations about
a given axis $\sigma$. We identify $so(3)$ with $\R^3$ using the cross product
and take the standard Euclidean inner product as our inner product on $\R^3$.
The horizontal subspaces are one dimensional at points without continuous 
isotropy, while all points with continuous isotropy are elements of the 
normalizer of $H$. Hence the curvature is identically zero. 

If \ $\hat {}: \R^3 \to so(3)$ is the standard identification of a 
three--vector with a skew--symmetric matrix, then 
$\mu(\hat \xi g) = (\la \sigma, \xi \ra, - \la g \sigma, \xi \ra)$, with
inertia factor
\[
\chi(g) = \lp \begin{array}{cc}
1 & - r(g) \\ - r(g) & 1 \end{array} \rp,
\]
where the invariant function $r: SO(3) \to \R$ is given by 
$r(g) := \la \sigma, g \sigma \ra$. The vectors $\nu_\pm := (1, \pm 1)$ are 
eigenvectors of $\chi(g)$, with eigenvalues $\lambda_\pm(g) = 1 \mp r(g)$. 
$g_\pm \in SO(3)$ has nontrivial isotropy if $g_\pm \sigma = \pm \sigma$, and
hence $\fg_{g_\pm} = \mbox{span}\{\nu_\pm\}$ and 
$\chi(g_\pm) = \nu_\mp \nu_\mp^T$. If we take $\pi = \half \chi(g_\pm)$
as our projection, then $\pi \, \chi(g) = \lambda_\mp(g) \pi$. Thus we can 
take $\iota(g) = \frac 1 {\lambda_\mp(g)} \Id$ as our restricted 
pseudo--inverse of $\chi$ and Corollary \ref{Abel_curv} implies there is
a slice $S$ through $g_\pm$ with 
\[
T_g S = \Xi_g 
= \lcb \hat \xi g : \la \xi, \sigma \pm g \sigma \ra = 0 \rcb.
\]

We can explicitly construct the slice $S$, using the Cayley transform
\[
\cay{\eta} = (\Id - \hat \eta/2)^{-1}(\Id + \hat \eta/2)
\]
from $\R^3 \approx so(3)$ to $SO(3)$ to define a coordinate map  
$\psi: \eta \mapsto \cay{\eta} g_0$ taking an open ball $B_r$ about 0 in $\R^3$ 
onto an open neighborhood of $g_0$ in $SO(3)$. A curve $\eta(s) \in \R^3$ with 
$\eta(0) = 0$ and $\norm{\eta(s)} < r$ for all $s$ lies in $S$ if 
$T_{\eta(s)} \psi \eta'(s) \in \Xi_{\eta(s)}$ for all $s$. This holds if and
only if
\beq{inter_calc}
0 = \tilde \nu(\psi(\eta(s))) T_{\eta(s)} \psi \eta'(s) 
= \la \sigma \pm \psi(\eta(s)) \sigma, \triv \psi(\eta(s)) \eta'(s) \ra.
\eeq
To further simplify this expression, we make use of the identities
$\Id + \cay{\eta} = 2 (\Id - \half \hat \eta)^{-1}$, and hence
\[
\sigma \pm \psi(\eta) \sigma = (\Id + \cay{\eta}) \sigma 
= (\Id - \half \hat \eta)^{-1} \sigma,
\]
and 
$\triv_\eta \mbox{cay} = \frac 1 {1 + \norm{\eta/2}^2} (\Id + \half \hat \eta)$.
Inserting these expressions into \for{inter_calc}, letting $\eta = \eta(s)$,
yields the condition
\[
0 = \la (\Id - \half \hat \eta)^{-1} \sigma, 
	(\Id + \half \hat \eta) \eta'(s) \ra
= \la (\Id + \half \hat \eta)^T (\Id - \half \hat \eta)^{-1} \sigma, 
	\eta'(s) \ra
= \la \sigma, \eta'(s) \ra.
\]
Thus $S = \psi(B_r \cap \sigma^\perp)$ is a slice through $g_0$.
\exampleend

\section{$\beta$--relative equivariant partial connections}
\label{beta_equiv}

We now relax the equivariance condition on the horizontal projections to 
allow for isotropy, requiring equivariance only with respect to some elements 
of the group. This is, in fact, the setting in which we originally introduced
partial connections; see \cite{LN00, LO, LN02}. Although the details of 
some calculations are more complicated in the $\beta$--equivariant setting,
the central constructions and underlying strategies are essentially identical
to those used in the fully equivariant setting. 

We shall say that a map $\beta: G \times M \to G$ is a 
{\em slip map} if $\beta(g, m) \cdot m = g \cdot m$ for all $g \in G$ and 
$m \in M$. If $G$ acts on manifolds $M$ and $N$ and $\beta$ is a slip map 
for the action on $M$, we shall say that a map $F: M \to N$ is 
{\em $\beta$--relative equivariant} if 
\beq{beta_eq}
F(\beta(g, m) \cdot m) = \beta(g, m) \cdot F(m)
\eeq
for all $g \in G$ and $m \in M$. 
Note that \for{beta_eq} can also be
expressed in the form $F(g \cdot m) = \beta(g, m) \cdot F(m)$. 
(See \cite{O99} for a more general treatment of relative equivariance.)

The constructions of partial connections, dual and partial connection forms, 
etc. all carry over to the $\beta$--equivariant setting. 
Analogs of Propositions 1--3 hold for $\beta$--equivariant
partial connections and forms. The inertia factor of a $\beta$--equivariant
dual connection form is $\beta$--equivariant. The definition of the curvature 
of a $\beta$--equivariant connection form or dual connection form is 
entirely analogous to that of an equivariant form. Note, however, that 
Lemma \ref{good_chi} need not hold for $\beta$--equivariant dual connection 
forms with nontrivial slip maps. Expressions analogous to the structure 
equations of Proposition \ref{structure} can be obtained by means of 
straightforward, but rather tedious, calculations. Development of results
analogous to those of \S \ref{almost_hor} for relatively equivariant 
connections will be the subject of future work.

To motivate the introduction of the notion of $\beta$--equivariance, 
we first relate trivial principal bundles to moving frames, then argue 
that $\beta$--equivariance, rather than full equivariance, is the most 
that can be expected in the setting of a natural extension of moving
frames to general actions. (See \cite{LN00, LO, LN02}.)

Given a Lie group $G$ acting on a manifold $M$,
an equivariant map $\rho: M \to G$, i.e. a map satisfying
\[
\rho(g \cdot m) = g \rho(m) \qquad \mbox{(left moving frame) \qquad or \qquad}
\rho(m) g \qquad \mbox{(right moving frame),}
\]
is called a {\em moving frame}. It is known that the existence of a 
moving frame implies that the group action is (locally) free. 
(See \cite{FOI, FOII}.) Given a
global cross--section $\cS$ of a manifold $P$ with a free $G$ action, i.e.
a transverse submanifold $\cS$ such that for each $p \in P$ there is a unique
$g \in G$ and $\tilde p \in \cS$ such that $p = g \cdot \tilde p$ if 
and only if there is a moving frame $\rho$ on $P$. 
If the action of $G$ on $M$ is globally free, then the existence of
a moving frame implies that $M$ is a trivial 
principal bundle. The base manifold can be naturally identified with 
the associated global cross--section $\sigma(M)$. It follows that in this
situation, if $\rho: M \to G$ is a moving frame, then the right 
trivialization $\triv_m \rho = d_m(R_{\rho(m)} \circ \rho)$ of the 
linearization of $\rho$ is a flat connection form on the principal bundle $M$. 
As we shall discuss below, the trivialized linearization of a generalization of
moving frames to non--free actions yields a relatively equivariant 
partial connection form.

\exam{A moving frame and connection on $US^2$}

The rotation group $SO(3)$ acts transitively on $S^2$ and freely and 
transitively on the unit tangent bundle $US^2 = \lcb u \in T S^2 : 
\norm{u} = 1 \rcb$. We consider the map 
\beq{sphere_mf}
\rho(u) := (m, u, m \times u)
\eeq
taking $u \in U_m S^2$
to the orthogonal matrix determined by the positively oriented 
orthonormal basis $\lcb m, u, m \times u \rcb$ is a moving frame. (Note 
that by an abuse of notation, we will regard $u$ both as a tangent vector
to the sphere at $m$ and as a unit vector in $\R^3$.) 

Since $\rho(g \, u) = (g \, m, g \, u, g (m \times u)) = g \, \rho(u)$
for any $g \in SO(3)$, $\rho$ is a left moving frame. 

The trivialized linearization $\triv \rho$ of $\rho$ satisfies
\beq{conn_unit_sphere}
\triv \rho(\du) = m \times \dm + \la u \times \du, m \ra m,
\eeq
where $\du \in T_u US^2$, with $m = \pi(u)$ and $\dm = d \pi \, \du$.
(Here $\pi: US^2 \to S^2$ denotes the canonical projection.) This can
be seen as follows: Let $u(\e)$ be a curve in $US^2$
through $u$ tangent to $\du$ and set $m_\e = \pi(u(\e))$. If we  
identify each of these points with a vector in $\R^3$, then 
differentiating the relations $\norm{m(\e)} = \norm{u(\e)} = 1$ and
$\la m(\e), u(\e) \ra = 0$ yields
\[
\la m, \dm \ra = \la u, \du \ra = \la \dm, u \ra + \la m, \du \ra = 0.
\]
Setting $z = m \times \dm + \la u \times \du, m \ra m$, we obtain
\[
d \rho(\du) 
= (\dm, \du, \dm \times u + m \times \du)
= (z \times m, z \times u, z \times (m \times u))
= \hat z \rho(m),
\]
where $\hat z$ denotes the skew--symmetric matrix satisfying
$\hat z y = z \times y$ for any $y$, $z \in \R^3$. 
\rmkend

In \cite{LN00, LN02, LO} we introduced the notion of a partial moving frame,
which is equivariant modulo isotropy. 
To motivate this extension, we observe that full equivariance holds if and 
only if $\rho(g \cdot m) \rho(m)^{-1} = g$ for all $m \in M$ and $g \in G$
and relax this condition to allow for isotropy. Assume that $G$ acts on
$M$ on the left. (If $G$ acts on the right, then the roles of left and 
right should be reversed throughout.)
We shall say that a (smooth) map $\phi: M \to G$ is a (left) 
{\em partial moving frame} if for any $g \in G$ and $m \in M$ 
\beq{partial moving}
\phi_g(m) := \phi(g \cdot m) \phi(m)^{-1} = g \quad \mbox{mod} \quad G_m.
\eeq
The quotient is with respect to right cosets; thus
\for{partial moving} is equivalent to 
\[
\phi_g(m) \cdot m = g \cdot m
\]
for all $g \in G$ and $m \in M$. 
(A {\em \pff} \ on a submanifold $\cS$ of a manifold $M$ with a $G$ action
is a map $\phi: \cS \to G$ satisfying \for{partial moving} for any $m \in \cS$ 
and any $g \in G$ such that $g \cdot m \in \cS$.)

To further motivate this construction, let us relate partial moving frames to 
cross--sections. Each partial moving frame $\phi: M \to G$ determines a global 
cross--section as follows. Define the map 
$\pi_\phi: M \to M$ by $\pi_\phi(m) := \phi(m)^{-1} \cdot m$.
Condition \for{partial moving} states that for any $m \in M$ and $g \in G$, 
there exists $h \in G_m$ satisfying $\phi(g \cdot m) = g \, h \, \phi(m)$,
and thus
\[
\pi_\phi(g \cdot m) 
= (g \, h \, \phi(m))^{-1} \cdot (g \cdot m) 
= \pi_\phi(m). 
\]
Thus each orbit $G \cdot m$ intersects the image of $\pi_\phi$ at precisely one point.
However, a global section does not uniquely determine a partial moving frame, 
since any (smooth) map $\iota: M \to G$ satisfying 
$\iota(m) \in G_m$ for all $m \in M$ determines a new partial moving frame 
$\tilde \phi(m) := \phi(m)\iota(m)$ satisfying $\pi_{\tilde \phi} =
\pi_\phi$. 

Given a partial moving frame 
on a bundle $\pi: N \to M$ with an equivariant projection $\pi$, we can 
construct a partial moving frame on $M$ using a section.

\begin{Proposition}
Let $M$ and $N$ be manifolds with $G$ actions. If
\begin{itemize}
\item
$\pi: N \to M$ is equivariant 
\item
$\sigma: M \to N$ satisfies $\pi \circ \sigma = \mbox{\rm id}$
and $\sigma(G \cdot m) \subseteq G \cdot \sigma(m)$ for all $m \in M$, 
\end{itemize}
then a partial moving frame $\phi$ on $N$ determines a partial moving frame
$\tilde \phi := \phi \circ \sigma$ on $M$.
\end{Proposition}
\prf
The inclusion $\sigma(G \cdot m) \subseteq G \cdot \sigma(m)$ 
implies that $\sigma(g \cdot m) = k \, \sigma(m)$ for some $k \in G$. Thus
equivariance of $\pi$ yields
\[
g \cdot m = \pi(\sigma(g \cdot m)) 
= \pi(k \cdot \sigma(m)) 
= k \cdot \pi(\sigma(m)) 
= k \cdot m,
\]
and hence $g^{-1} k \in G_m$. Combining this with
\[
k^{-1} \tilde \phi_k(m) = k^{-1} \phi_k(\sigma(m))
\in G_{\sigma(m)} \subseteq G_m
\]
gives $g^{-1} \tilde \phi_g(m) \in G_m$.
\prfend

\exam{Partial moving frames and connection forms on subsets of $S^2$}
The projection $\pi$ from $US^2$ to $S^2$ is clearly 
equivariant. Hence any unit vector field $Y$ on a submanifold $M$ of $S^2$ 
determines a \pff \ $\phi = \rho \circ Y$ on $M$. 
It follows immediately from \for{conn_unit_sphere} that the partial moving 
connection form $\triv \phi$ determined by $\phi$ is
\beq{fuzzconn_sphere}
\triv \phi(\dm) = m \times \dm + \la Y(m) \times (dY(\dm)), m \ra m.
\eeq
The map $\phi_g$ associated to $g \in SO(3)$ is
$\phi_g(m) = g \, \exp(\theta(g, m) \, m)$, where
$\theta(g, m)$ denotes the angle between $g^{-1} Y(g \, m)$ and $Y(m)$. 
Thus $g^{-1} \phi_g(m)$ measures the failure of the vector field
$Y$ to commute with the action of $g$ at $m$. To see this,
let $u = g^{-1} Y(g \, m) \in U_m S^2$ and $\theta = \theta(g, m)$. Then 
$\phi(g \, m) = \rho(Y(g \, m)) = g \, \rho(u)$
and, using $\la m, u \ra = \la g \, m, Y(g \, m) \ra = 0$  and
$\la m \times u, m \times Y(m) \ra = \la u, Y(m) \ra$, we obtain
$\phi(m)^T \rho(u) = \exp(\theta \, e_1)$,
where $e_1$ denotes the first standard Euclidean basis vector. Hence
\[
\phi_g(m) 
= \phi(g \, m) \phi(m)^{-1} 
= g \, \rho(u) \phi(m)^{-1}.
\]
Note that if $m(t)$ is a unit speed curve in $S^2$, then 
\[
\triv \rho(\ddot m) = m \times \dot m + \la \dot m \times \ddot m, m \ra m
= m \times \dot m + k_g(m) m,
\]
where $k_g(m(t))$ denotes the geodesic curvature of $m(t)$.
Similarly, if the unit vector field $Y$ is the normalization of a nonzero
vector field $X$ on some set $M \subset S^2$,
i.e. $Y(m) = X(m)/\norm{X(m)}$, then \for{fuzzconn_sphere} yields
\[
\triv \phi(X(m)) 
= m \times X(m) + k_g(m) \,\norm{X(m)} m.
\qquad \rmkend
\]

The relevance of partial moving frames to partial connections is 
established in the following result.

\begin{Theorem}
\label{lin_pmf}
The trivialized linearization $\triv \phi$ of a partial moving frame 
$\phi: M \to G$ is a $\beta$--equivariant partial connection form, with 
associated slip map $\beta(g, m) = \phi_g(m)$. 

If we define $\pi_\phi(m) := \phi(m)^{-1} \cdot m$, then 
the partial connection $\Gamma$ determined by $\triv \phi$ satisfies
\[
\Gamma|_m = d (\Phi_{\phi(m)} \circ \pi_\phi)(T_m M)
\]
for all $m \in M$.
\end{Theorem}

We now develop some consequences of $\beta$--equivariance, some of which 
we will use in the proof of Theorem \ref{lin_pmf}.

\begin{Lemma}
\label{beta_var}
\ \vspace{-.1in} \\
\begin{enumerate}
\item
If $\beta_g: M \to G$ satisfies $\beta_g(m) \cdot m = g \cdot m$
for all $m \in M$, then
\beq{beta_eq1}
d_m \Phi_g - d_m \Phi_{\beta_g(m)}
	= d \hat \Phi_{g \cdot m} \circ \triv_m \beta_g
\eeq
and
\beq{local1}
\range{(\triv \beta_g \circ d_e \hat \Phi_m - \Ad g + \Ad {\beta_g(m)}}) 
\subseteq \fg_{g \cdot m}
\eeq
for all $m \in M$. 
\item 
A $\fg$--valued one--form $\alpha$ is $\beta$--equivariant modulo isotropy
if and only if 
\beqa
\Phi_g^*\alpha(v) 
&=& (\Ad {\beta(g, m)}\alpha + \triv \beta_g)(v)
	\quad \mod \quad \fg_{g \cdot m} \\
&=& (\Ad g \alpha + \triv \beta_g \circ (\Id - \proj \alpha)(v)
	\quad \mod \quad \fg_{g \cdot m},
\eeqa
for all $g \in G$, $m \in M$, and $v \in T_M M$.
In particular, $\alpha$ is fully equivariant modulo isotropy if and only if 
$\triv \beta_g \circ (\Id - \proj \alpha)$ maps $T_m M$ into $\fg_{g \cdot m}$ 
for all $g \in G$ and $m \in M$. 
\item
A $\fg^*$--valued one--form $\mu$ is $\beta$--equivariant if and only if
\[
\Phi_g^*\mu(v) = (\Ads {\beta(g, m)^{-1}} \mu
	+ \chi(g \cdot m) \triv \beta_g)(v) 
\]
for all $g \in G$, $m \in M$, and $v \in T_mM$.
\end{enumerate}
\end{Lemma}
\prf
(i) Given $g \in G$, $m \in M$, and $v \in T_m M$, set $h = \beta_g(m)$ 
and $\zeta = \triv \beta_g(v)$. To prove (i), note that
linearizing $\beta_g(m) \cdot m = g \cdot m$ with respect to $m$ yields
\[
d \Phi_g(v) = d_h \hat \Phi_m (d\beta_g(v) + d\Phi_h(v)).
\]
Hence
\[
d \hat \Phi_m (d \beta_g(v)) 
= d (\hat \Phi_m \circ R_h) \zeta 
= d \hat \Phi_{h \cdot m} \zeta 
= \zeta_M(g \cdot m).
\]
Regrouping terms yields \for{beta_eq1}.

Combining \for{beta_eq1} and the identity 
$d_e(\Phi_g \circ \widehat \Phi_m) = d \widehat \Phi_{g \cdot m} 
	\circ \Ad g$, for all $g \in G$ and $m \in M$, yields
\[
d \widehat \Phi_{g \cdot m} \circ \triv \beta_g 
	\circ d_e \widehat \Phi_m
= (d \Phi_g - d \Phi_h) \circ d_e \widehat \Phi_m
= d \widehat \Phi_{g \cdot m} \circ (\Ad g - \Ad h),
\]
and hence $\triv \beta_g \circ d_e \widehat \Phi_m = (\Ad g - \Ad h)$ modulo 
$\fg_{h \cdot m}$. 
\medskip

\noindent
(ii) and (iii)
If $\alpha$ is $\beta$--equivariant modulo isotropy, then
\beqan
\Phi_g^* \alpha(v) &=& \alpha(d \Phi_g (v)) \nonumber \\
&=& \alpha(d \Phi_{\beta(g, m)}(v) + d \actM {g \cdot m} \triv \beta_g(v)) 
	\label{local2} \\
&=& \Ad {\beta(g, m)} \alpha(v) + \triv \beta_g(v)
	\quad \mod \fg_{g \cdot m} \nonumber  
\eeqan
yields the first equality of (ii). Analogously, if $\mu$ is a 
$\beta$--equivariant $\fg^*$--valued form, then 
\[
\Phi_g^* \mu(v) 
= \mu(d \Phi_{\beta(g, m)}(v) + d \actM {g \cdot m} \triv \beta_g(v)) 
= \Ads {\beta(g, m)} \mu(v) + \chi(m) \triv \beta_g(v)).
\]
To prove the second equality of (ii), we combine \for{local1} and \for{local2},
obtaining
\[
\Phi_g^* \alpha(v) 
= \Ad g \alpha(v) + \triv \beta_g((\Id - d \actM m \circ \alpha)(v))
	\quad \mod \fg_{g \cdot m}.
\]
Since $\proj \alpha|T_m M = d \actM m \circ \alpha_m$, this completes
the proof.
\prfend

\prfof{Theorem \ref{lin_pmf}}
Condition \for{partial moving} can be written as 
$R_{\phi(m)^{-1}}(\phi(\widehat \Phi_m(g))) = g$ modulo $G_m$. 
Linearizing with respect to the group element yields 
\[
\triv \phi(\xi_M(m)) 
= d (R_{\phi(m)^{-1}} \circ \phi \circ \widehat \Phi_m)\xi
= \xi \quad \mod \fg_m
\]
for all $\xi \in \fg$ and $m \in M$. 

We now show that the trivialized linearization $\triv \phi$ of a 
partial moving frame $\phi: M \to G$ satisfies
\[
\triv \phi(d \Phi_g (v)) 
= \Ad {\phi_g(m)} \triv \phi (v) + \triv \phi_g(v)
\]
for all $v \in T_m M$,
and hence, applying Lemma \ref{beta_var}, is $\beta$--equivariant.
Given $v \in T_m M$, let $m(\e)$ be a curve through $m$ in $M$
satisfying $\ddeps m(\e) |_{\e = 0} = v$. Then
\beqa
\triv \phi(d \Phi_g(v)) 
&=& d(R_{\phi(g \cdot m)^{-1}} \circ \phi \circ \Phi_g)(v) \\
&=& \frac {d \ } {d \e} \phi(g \cdot m_\e) \phi(g \cdot m)^{-1}
	|_{\e = 0} \\
&=& \frac {d \ } {d \e} \phi_g(m(\e)) \phi(m(\e)) \phi(m)^{-1}
	 \phi_g(m)^{-1}|_{\e = 0} \\
&=& \Ad {\phi_g(m)} \triv \phi (v) + \triv \phi_g(v).
\eeqa
The relatively equivariant analog of Proposition \ref{conn_relations1}
thus implies that $\triv \phi$ determines a $\beta$--equivariant
projection onto the tangent spaces of the group orbits.

Linearizing the relation $\phi(m) \cdot \pi_\phi(m) = m$ yields 
\[
d \Phi_{\phi(m)} \circ d_m \pi_\phi + d \actM {\pi_\phi(m)} \triv \phi = \Id,
\] 
and hence
\[
d_m \pi_\phi 
= d \Phi_{\phi(m)^{-1}} \circ (\Id - d \actM {\pi_\phi(m)} \circ \triv \phi)
= d \Phi_{\phi(m)^{-1}} \circ \proj \Gamma,
\] 
where $\Gamma$ is the partial connection determined by $\triv \phi$. Thus
$\range{d_m \pi_\phi} =  d \Phi_{\phi(m)^{-1}}(\Gamma|_m)$. 
\prfend

\medskip

\noindent 
{\bf Acknowledgements.}
The authors would like to thank Viktor Ginzburg and Hans Munthe-Kaas for 
many helpful discussions and suggestions.

\bibliographystyle{unsrt}
\bibliography{fuzzy}

\section*{Appendix}

\restate{Theorem \ref{conn_relations4}}{
\begin{enumerate}
\item
A dual connection form $\mu$ determines a partial connection
$\parconn = \kernel{\mu}$ and inertia factor
$\chi = \mu \circ d_e \actM {}$.

\item
An equivariant (singular) differential system $\parconn$ satisfying
$T M = \tfg \oplus \parconn$ is a partial
connection if there is inertia factor $\chi$ such that the equivariant
$\fg^*$--valued one--form $\mu$ given by
\[
\mu|{\parconn} := 0
\sands
\mu \circ d_e \actM m := \chi(m)
\quad \mbox{for all $m \in M$}
\]
is smooth, and hence a dual connection form.

\item
A $\fg$--valued one--form $\alpha$ is a partial connection form if
there is an inertia factor $\chi$ such that $\mu = \chi \, \alpha$ is a dual
connection form with inertia factor $\chi$.
\end{enumerate}
}
\prf
(i) Equivariance of $\mu$ implies that $\Gamma$ is an equivariant differential system;
hence $\Gamma$ is a partial connection.
\medskip

\noindent
(ii) $\kernel{(\mu \circ T_e \actM m)} = \kernel{\chi(m)} = \fg_m$ for all
$m \in M$ implies that $\mu|\tfg$ is one-to-one, and hence
$\kernel{\mu} = \Gamma$. Hence $\Gamma$ is a partial connection.
Any tangent vector $v \in T_m M$ satisfies $v = \xi_M(m) + u$ for some
$\xi \in \fg$ and $u \in \parconn|_m$. Equivariance of $\parconn$ implies
that $d \Phi_g u \in \parconn|_{g \cdot m}$, while
\[
d_e (\Phi_g \circ \actM m)
= d \actM {g \cdot m} \circ \Ad g
= \chi(g \cdot m) \circ \Ads g
= \Ads {g^{-1}} \chi(m)
\]
mplies that
\[
\Phi^*_g \mu(v)
= \mu(d \Phi_g (\xi_M(m) + u))
= \Ads {g^{-1}} \chi(m) \xi
= \Ads {g^{-1}} \mu(v)
\]
for any $g \in G$. Thus $\mu$ is equivariant and, hence, if smooth, is a dual
connection form.
\medskip

\noindent
(iii) If $\mu = \chi \, \alpha$
is a dual connection form with associated inertia factor $\chi$, then
\[
\chi(m)  = \mu \circ d_e \actM m
= \chi(m) \circ \alpha \circ d_e \actM m
\]
and $\kernel{\chi(m)} = \fg_m$ imply that $\alpha$ satisfies \for{proj_range}.
Equivariance of $\mu$ and $\chi$ imply that
\[
\chi \, \Ad {g^{-1}} \Phi_g^* \alpha
= \Ads g \Phi_g^* \mu
= \mu
= \chi \, \alpha
\]
and thus
\[
\range{(\Ad {g^{-1}} \Phi_g^* \alpha - \alpha)_m}
        \subseteq \kernel{\chi(m)} = \fg_m,
\]
i.e. that $\alpha$ is equivariant modulo isotropy. Hence Proposition
\ref{conn_relations1} implies that $\alpha$ is a partial connection form.
\prfend

\restate{Lemma \ref{slice_cond}}{
\begin{enumerate}
\item
Given an involutive $\gmo$--equivariant differential system
$\Delta \supseteq \widetilde {\fgmo}$ on a $\gmo$--invariant neighborhood
of $m_0$, there
is an $\gmo$--invariant integral submanifold of $\Delta$ containing $m_0$.
\item
If $G$ acts properly and $S$ is a $\gmo$--invariant submanifold of $M$ containing
$m_0$ and satisfying
$T_{m_0} M = T_{m_0} S \oplus \tfg|_{m_0}$, there is a neighborhood
$S_0$ of $m_0$ in $S$ such that $S_0$ is a slice through $m_0$.
\end{enumerate}
}
\prf
\noindent (i)
Let $N$ be a neighborhood of $m_0$ in the maximal integral manifold of 
$\Delta$ containing $m_0$ and let $v \in T_{g \cdot m}(G \cdot N)$ for
some $g \in G$ and $m \in N$. There are curves $g(\epsilon) \in G$ and 
$m(\epsilon) \in N$, with $g(0) = g$, $m(0) = m$, such that
\[
v = \frac {d \ }{d \epsilon} g(\epsilon) \cdot m(\epsilon)|_{\epsilon = 0}
= \xi_M(g \cdot m) + d_m \Phi_g(w),
\] 
where $\xi = \frac {d \ }{d \epsilon} g(\epsilon)|_{\epsilon = 0}$ and
$w = \frac {d \ }{d \epsilon} m(\epsilon)|_{\epsilon = 0} \in \Delta|_m$. 
Hence equivariance of $\Delta$ implies that
\[
v = d_m \Phi_g((\Ad {g^{-1}} \xi)_M(m) + w) 
\in d_m \Phi_g(\Delta|_m) 
= \Delta|_{g \cdot m}.
\]
Thus $T_{g \cdot m}(G \cdot N) \subseteq \Delta|_{g \cdot m}$. Maximality
of $N$ implies that $G \cdot N = N$.
\medskip

\noindent (ii)
To show that $S$ is a slice, it suffices to show that there is some 
neighborhood $S_0$ of $m_0$ in $S$ such that $(g \cdot S_0) \cap S_0 
\neq \emptyset$ implies $g \in \gmo$. Our proof largely follows that 
given in Duistermaat and Kolk. Consider a sequence $(g_j, m_j)$ in $G \times S$ such that
$g_j \cdot m_j \in S$ for all $j$ and 
\[
\lim_{j \to \infty} g_j \cdot m_j = m_0 = \lim_{j \to \infty} m_j.
\]
Passing to a convergent subsequence if necessary, let 
$g = \lim_{j \to \infty} g_j$; $g \cdot m_0 = \lim_{j \to \infty} g_j \cdot m_j 
= m_0$ implies that $g \in \gmo$. If we set $h_j := g^{-1} g_j \not \in \gmo$, 
so that $\lim_{j \to \infty} h_j = e$, then the $\gmo$--invariance of $S$ 
implies that $h_j \cdot m_j \in S$ for all $j$. 

The decompositions
\[
\fg = \fg_{m_0} \oplus \fh
\sands
T_{m_0} M = \tfg|_{m_0} \oplus \Xi|_{m_0}
= \tilde {\fh}|_{m_0} \oplus T_{m_0} S
\]
and the Inverse Function Theorem imply that there are 
$\gmo$--invariant neighborhoods $\cV$ of 0 in $\fh$, $\cW$ of $e$ in $\gmo$, 
and $S_0$ of $m_0$ in $S$ such that the equivariant maps
$\Psi_G: \cV \times \cW \to G$ and $\Psi_M: \cV \times S_0 \to M$ given by
\[
\Psi_G(\eta, g) := \exp(\eta) g
\sands
\Psi_M(\eta, m) := \exp(\eta) \cdot m
\]
are diffeomorphisms onto their images. 
For any $\eta \in \cV$, $g \in \cW$, and $m \in S_0$, we have
\[
\Psi_G(\eta, g) \cdot m = \Psi_M(\eta, g \cdot m).
\]

For sufficiently large $j$, we have $m_j \in S_0$, $h_j \cdot m_j \in S_0$, and 
$h_j = \Psi_G(\eta_j, k_j)$ for some $k_j \in \cW$ and $\eta_j \in \cV$. Hence
\[
\Psi_M(0, h_j \cdot m_j) = h_j \cdot m_j 
= \Psi_G(\eta_j, k_j) \cdot m_j 
= \Psi_M(\eta_j, k_j \cdot m_j).
\]
Injectivity of $\Psi_M$ implies that $\eta_j = 0$ and hence $h_j = k_j \in \gmo$.
\prfend

Before proving Proposition \ref{adapt_exist}, we establish a close relationship 
between a certain class of nondegenerate adaptors and slices. Given an
adaptor $\phi$, if the maps $\bleh m := \triv_e (\phi \circ \actM {m})$ satisfy
$\range{(\Id - \bleh m)} \subseteq \fgmo$ for all $m \in \phi^{-1}(\gmo)$,
then $\phi$ is {\em transversal}. Here $\triv f$ denotes the right 
trivialization of the linearization of a map $f: M \to G$, i.e.
$\triv_m f = d_m (R_{f(m)^{-1}} \circ f)$. The following proposition shows that 
if the action is proper, then slices determine transversal adaptors and vice 
versa.

\begin{Proposition}
\label{frame_stuff}
Assume that $G$ acts properly on $M$.
\begin{itemize}
\item
If $S$ is a slice through $m_0$, there is a transversal adaptor $\phi$ 
for $m_0$ such that $\phi^{-1}(\gmo)$ is a neighborhood of $m_0$ in $S$. Given 
any $h \in \gmo$ and any $\gmo$--invariant complement $\fh$ of $\fg_{m_0}$, 
$\phi$ can be chosen so that $\phi(m_0) = h$ and $\range{\bleh {m_0}} = \fh$. 
\item
If $\phi$ is a transversal adaptor for $m_0$, then some neighborhood of 
$m_0$ in $\phi^{-1}(\gmo)$ is a slice through $m_0$ and
$\range{\bleh {m_0}}$ is a $\gmo$--invariant complement of $\fg_{m_0}$.
\end{itemize}
\end{Proposition}

\prf
\noindent (i) 
Transversality of $\phi$ implies that
\[ 
\range{d \phi_m} + T_{\phi(m)} \gmo
\supseteq d R_{\phi(m)} (\range{\bleh m} + \fg_{m_0})
= d R_{\phi(m)} \fg 
= T_{\phi(m)} G
\]
for all $m \in S := \phi^{-1}(\gmo)$. Hence $S$
is a submanifold of $M$; $\gmo$--equivariance of $\phi$ implies that $S$ is
$\gmo$--invariant. Since $\phi$ is transversal, the restriction of
$\bleh {m_0}$ to $\fh := \range{\bleh {m_0}}$ is an isomorphism, and hence
\[
T_{m_0} M = \kernel{\triv_{m_0} \phi} \oplus \tilde \fh|_{m_0}
= T_{m_0} S \oplus \tfg|_{m_0}.
\]
Thus Lemma \ref{slice_cond} implies that some neighborhood $S_0$ of $m_0$ in
$S$ is a slice through $m_0$. If $m \in S$, then $\phi(m) \in \gmo$, and hence
$G_{\phi(m)^{-1} \cdot m} \subseteq \gmo$ implies that $G_m \subseteq \gmo$.
Finally, $\gmo$--equivariance of $\phi$ and $\actM {m_0}$ imply that $\fh$
is $\gmo$--invariant.

\noindent (ii) 
Restricting $S_0$ if necessary, there are $\gmo$--invariant neighborhoods $\cV$
of 0 in $\fh$, $\cW$ of $e$ in $\gmo$, $\cG$ of $e$ in $G$, and $\cU$ of $m_0$
in $M$ such that the equivariant maps
$\Psi_G: \cV \times \cW \to \cG$ and $\Psi_M: \cV \times S_0 \to \cU$ given by
$\Psi_G(\eta, h) := \exp(\eta) h$ and $\Psi_M(\eta, m) := \exp(\eta) \cdot m$
are diffeomorphisms satisfying
\[
\Psi_G(\eta, h) \cdot m = \Psi_M(\eta,  h \cdot m)
\]
for all $\eta \in \cV$, $h \in \cW$, and $m \in \cU$. 
If we set $\phi(\Psi_M(\eta, m)) := \exp(\eta)$, then $\phi|{S_0} \equiv e$ and
\[
\phi(h \cdot \Psi_M(\eta, m))
= \phi(\Psi_M(\Ad h \eta, h \cdot m))
= \exp(\Ad h \eta)
= h \, \exp(\eta) \, h^{-1}
\]
for $h \in \gmo$. If $m \in \cU$ and $g = \Psi_G(\eta, h) \in \cG$, then
\[
g^{-1} \cdot \phi(g \cdot m)
= (\exp(\eta) h)^{-1} \phi(\Psi_M(\eta, h \cdot m))
= h^{-1} \in \gmo.
\] 
Hence we can define $\psi_m: \cG \to \gmo$ by 
$\psi_m(g) := g^{-1} \phi(g \cdot m)$, with
\[
\triv_e \psi_m = \Ad {\phi(m)}(\Id - \bleh m).
\]
Thus $\phi(m) \in \gmo$ implies that $\range{(\Id - \bleh m)}
\subseteq \fgmo$, and hence $\phi$ is transversal.
\prfend 

\restate{Proposition \ref{adapt_exist}}
{
If $G$ acts properly, then there is an adaptor for any point in $M$. Given a
partial connection $\Gamma$ and a point $m_0 \in M$, there is an adaptor
$\phi$ for $m_0$ satisfying $\phi(m_0) = e$ and
$d \phi(\Gamma|_{m_0}) \subseteq \fgmo$.
}
\medskip

\prf
Equip $M$ with a $\gmo$--invariant metric and let $S'$ denote the image of
a $\gmo$--invariant neighborhood of $0$ in $\Gamma|_{m_0}$ under the 
associated exponential map. $S'$ is $\gmo$--invariant and satisfies 
\[
T_{m_0} S' \oplus \tfg|_{m_0} 
= \Gamma|_{m_0} \oplus \tfg|_{m_0} 
= T_{m_0} M.
\]
Hence Lemma \ref{slice_cond}.ii implies that some neighborhood $S$ of $m_0$ in $S'$
is a slice. Proposition \ref{frame_stuff} now implies that there is a
transversal adaptor $\phi$ for $m_0$ such that $\phi(m_0) = e$ and 
$\phi$ maps some neighborhood of $m_0$ in $S$ into $\gmo$; hence
\[
\trivl \phi(\Gamma|_{m_0}) = d \phi(T_{m_0}S) 
\subseteq T_{\phi(m_0)} \gmo = \fgmo.
\qquad \prfend
\] 

\restate{Corollary \ref{flat_reg_lcs}}{
If $G$ acts properly on $M$,
$\Omega$ equals zero near $m_0$, and there is an adaptor $\phi$ satisfying
\beq{full2}
G_m = \phi(m) \gmo \phi(m)^{-1} \qquad \mbox{for all $m$ near $m_0$},
\eeq
then $\Gamma$ is tangent to a local cross section through $m_0$.
}
\medskip

\prf
Assume that $\Gamma$ arises from a regular foliation.
Let $S'$ denote the leaf containing $m_0$. The decomposition
\[
T_{m_0} M = \tfg|_{m_0} \oplus \Gamma|_{m_0}
= \tilde {\fh}|_{m_0} \oplus T_{m_0} S'
\]
and the Inverse Function Theorem imply that there are neighborhoods
$\cV$ of 0 in $\fh$ and $S''$ of $m_0$ in $S'$ such that the map
$\Psi_M: \cV \times S'' \to M$ given by
\[
\Psi_M(\eta, m) := \exp(\eta) \cdot m
\]
is a diffeomorphism onto its image. Hence $G \cdot S'' \supseteq
\Psi_M(\cV \times S'')$ contains a neighborhood of $m_0$ in $M$.

To show that some neighborhood $S$ of $m_0$ in $S''$ is a local cross section
if there is an adaptor $\phi$ satisfying \for{full2},
it remains to be shown that $m \in S$ and $g \cdot m \in S$ implies
$g \in G_m$. The map $F : \cV \times \gmo \times S' \to G$ given by
\[
F(\eta, h, m) := \exp(\eta) \phi(m) h \phi(m)^{-1}
\]
satisfies $F(\eta, h, m) \cdot m = \Psi_M(\eta, m)$ and
\[
d_{(0, h, m_0)} F(\xi, d R_h \zeta, 0)
= d R_{\phi(m_0) h \phi(m_0)^{-1}}(\xi + \Ad {\phi(m_0)} \zeta)
\]
for any $h \in \gmo$. The decomposition $\fg = \fg_{m_0} \oplus \fh$
and the Implicit Function Theorem imply that there are neighborhoods
$\cG$ of $h_0$ in $G$ and $S$ of $m_0$ in $S''$ such that
the maps $\eta: \cG \times S \to \cV$ and $h: \cG \times S \to \gmo$
satisfy $F(\eta(g, m), h(g, m), m) = g$ for all $g \in \cG$ and $m \in S$.

Consider a sequence $(g_j, m_j)$ in $G \times S$
such that $g_j \cdot m_j \in S$ for all $j$ and
\[
\lim_{j \to \infty} g_j \cdot m_j = m_0 = \lim_{j \to \infty} m_j.
\]
Passing to a convergent subsequence if necessary, let
$g = \lim_{j \to \infty} g_j$; $g \cdot m_0 = \lim_{j \to \infty} g_j \cdot m_j
= m_0$ implies that $g \in \gmo$.
For sufficiently large $j$, we have $m_j \in S$, $g_j \cdot m_j \in S$, and
$g_j \in \cG$. Let $\eta_j = \eta(m_j, g_j)$ and $h_j = h(m_j, g_j)$. Then
\[
\Psi_M(0, g_j \cdot m_j) = F(\eta_j, h_j, m_j) \cdot m_j
= \Psi_M(\eta_j, m_j).
\]
Injectivity of $\Psi_M$ implies that $\eta_j = 0$ and $g_j \in G_{m_j}$.
\prfend

\end{document}